\documentclass[journal]{IEEEtran}
\usepackage{amsfonts}
\usepackage{amsmath}
\usepackage{amssymb}
\usepackage{amsthm}
\usepackage{booktabs}
\usepackage{multirow}
\usepackage{graphicx}
\usepackage{float}
\usepackage{overpic}
\usepackage{lineno,hyperref}
\usepackage{siunitx}
\usepackage{amstext}
\usepackage{longtable}
\usepackage{tabularx}
\usepackage{array}
\usepackage{enumitem}
\usepackage[linesnumbered,ruled]{algorithm2e}
\theoremstyle{definition}
\newtheorem{definition}{\textbf{{Definition}}}

\ifCLASSINFOpdf
\else
\fi
\begin{document}
	
\title{Hierarchical Prior Regularized Matrix \\ Factorization for Image Completion}

\author{Liyu Su~{~}% <-this % stops a space
\thanks{L. Su is with the Faculty
of Electronic and Information Engineering, Xi'an Jiaotong University, Xi'an,
710049 China. E-mail: luc.su@stu.xjtu.edu.cn. Copyright of this work may be transferred without notice, after which, the current version may no longer be accessible.}}% <-this % stops a space


\markboth{Journal of \LaTeX\ Class Files,~Vol.~14, No.~8, August~2015}%
{Shell \MakeLowercase{\textit{et al.}}: Bare Demo of IEEEtran.cls for IEEE Journals}



\maketitle

\begin{abstract}
The recent low-rank prior based models solve the tensor completion problem efficiently. However, these models fail to exploit the local patterns of tensors, which compromises the performance of tensor completion. In this paper, we propose a novel hierarchical prior regularized matrix factorization model for tensor completion. This model hierarchically incorporates the low-rank prior, total variation prior, and sparse coding prior into a matrix factorization, simultaneously characterizing both the global low-rank property and the local smoothness of tensors. For solving the proposed model, we use the alternating direction method of multipliers to establish our algorithm. Besides, the complexity and convergence are investigated to further validate the algorithm effectiveness. The proposed scheme is then evaluated through various data sets. Experiment results verify that, the proposed method outperforms several state-of-the-art approaches.
\end{abstract}


% \fi
\IEEEpeerreviewmaketitle

\section{Introduction}

\IEEEPARstart{T}{ensors}, as multidimensional arrays, attract intensive attention for their natural strength to represent data with high order structure, e.g., a color image is a third-order tensor characterized by two spatial modes, and one color mode respectively. For decades, tensors have emerged in image processing \cite{ip1}-\cite{ip4}, computer vision \cite{cv1}-\cite{cv3}, data mining \cite{dm1}-\cite{dm3}, and machine learning \cite{ml1}-\cite{ml4}, etc. However, due to the transmission and storage restrictions, incomplete tensors are common concerns in practice. To tackle the concerns, tensor completion methods are performed for estimating the missing entries of incomplete tensors.

Inspired by low-rank matrix completion schemes \cite{mc1}-\cite{mc4}, low-rank tensor completion methods are popular tools to exploit the global structures of incomplete tensors. Unlike the matrix rank, the definition of tensor rank is still an open problem. Recently, the tensor CANDECOMP/PARAFAC (CP) rank \cite{cp1}-\cite{cp4}, tensor Tucker rank \cite{tu1}-\cite{tu3}, tensor tubal rank \cite{tubal1}-\cite{tubal3}, and tensor train (TT) rank \cite{tt1}-\cite{tt3} are incorporated as the low-rank priors for tensor completion. In detail, the tensor CP rank is defined as the number of rank-one tensors. Unfortunately, minimizing the tensor CP rank is NP-hard yet \cite{cpnphard}. The tensor Tucker rank is defined as the ranks of factor matrices, which inevitably destroys the tensor internal structure. To alleviate this structure destruction, the tensor tubal rank uses Fourier transform to effectively characterize the tensor structure, but it is only available for third-order tensors. Furthermore, the tensor TT rank is proposed via the tensor train decomposition. Despite the improvement, the tensor TT rank lacks representation flexibility, since the border ranks have a fixed pattern.

In addition to the low-rank priors, the smoothness priors are imposed to describe the local patterns of incomplete tensors. Explicitly, the total variation (TV) is a typical prior to enhance piecewise smoothness \cite{tv1}-\cite{tv3}, minimizing the difference between neighboring tensor entries. Besides, the sparse coding (SC) is another efficient prior to induce the local smoothness \cite{sc1}-\cite{sc3}. Compared to the TV prior, the SC prior is often expressed through dictionaries, which are further categorized into  learned dictionaries and prespecified dictionaries. The learned dictionaries are mainly obtained via the dictionary learning, showing flexibility and adaptivity to the specific data. Nevertheless, the dictionary learning may be time-consuming when processing high order tensors. The prespecified dictionaries include various transform based dictionaries, e.g., contourlets, wavelets, discrete cosine transforms (DCT), etc. These dictionaries are simple to implement, and require limited computaional complexity.

Increasing researchers claim the exclusive employment of low-rank priors or smoothness priors may lead to compromised results for tensor completion. Concerning low-rank priors, each tensor rank describes the tensor structure from different perspectives. Joint utilization of low-rank priors may illustrate improvement \cite{cp+tu}. Furthermore, simultaneously using low-rank and smoothness priors to mine both global and local tensor features presents significant progress \cite{rank+smooth1}-\cite{rank+smooth3}. Although the existing approaches show promising development, they independently consider the low-rank and smoothness priors, which yields redundancy and increases complexity.

In this paper, we propose a novel model by simultaneously capturing the global low-rank property and the local smoothness to recover incomplete tensors. Consequently, the corresponding Lagrangian function is formulated, and we solve the function via the program of alternating direction method of multipliers (ADMM). To determine the subproblems in terms of the low-rank prior and the smoothness prior, we optimize one prior term while fix other prior terms at a time. Moreover, the convergence and complexity are also investigated to further depict the performance of the proposed algorithm. Experiments on different data sets show that, our proposed algorithm is superior to several state-of-the-art algorithms. In summary, the contributions of this paper are listed as follows:
\begin{enumerate}[leftmargin=*]
\item We establish a novel hierarchical prior regularized matrix factorization (HPMF) model for tensor recovery, which hierarchically adopts the low-rank prior, TV prior, and SC prior to explore both global structures and local patterns of tensors. This hierarchical prior framework can be generalized to different factorization schemes with other priors.
\item We propose an effective optimization algorithm to solve the HPMF model by the program of ADMM. Additionally, the complexity and convergence analyses are investigated, indicating our method shows acceptable computational efficiency with guaranteed convergence.
\item We utilize various data sets to test the performance of the proposed algorithm. The experimental results demonstrate that, our algorithm outperforms several state-of-the-art algorithms not only quantitatively but also visually.
\end{enumerate}

The remainder of this paper is organized as follows. Section \ref{sec-nota} introduces preliminaries involved in tensor completion. In section \ref{sec-model}, we demonstrate the proposed optimization model. The algorithm for the proposed model is determined in section \ref{sec-algo}. We conduct experiments in section \ref{sec-experi}. Finally, section \ref{sec-conclu} illuminates the conclusion.

\section{Preliminaries}
\label{sec-nota}
This section introduces notations, definitions and reviews related works for tensor completion methods.

\subsection{Notations}
In this article, we show scalars by lower case letters, e.g., $x$, $y$, $z$, and depict matrices as upper case letters, e.g., $X$, $Y$, $Z$. The $(i,j)$-th entry of matrix $X$ is represented via $X_{i,j}$. Tensors are denoted as calligraphic letters, e.g., $\mathcal{X}$, $\mathcal{Y}$, $\mathcal{Z}$. Precisely, an $N$th-order tensor is denoted as $\mathcal{X} \in \mathbb{R}^{I_{1} \times I_{2} \times \cdots \times I_{N}}$, where $I_n$, $n=1,2,\dots,N$, is the $n$th mode of $\mathcal{X}$. We demonstrate $\mathcal{X}_{i_1,i_2,\dots,i_N}$ for the $(i_1,i_2,\dots,i_N)$-th tensor entry of $\mathcal{X}$.

\subsection{Definitions}
The necessary definitions used in this paper are briefly provided as follows:

\begin{definition} [Frobenius norm]
The Frobenius norm of an $N$th order tensor $\mathcal{X}\in \mathbb{R}^{I_{1} \times I_{2} \times \cdots \times I_{N}}$ is defined as:
\begin{equation}
	\|\mathcal{X}\|_F=(\sum^{I_1}_{i_1=1}\sum^{I_2}_{i_2=1}\cdots\sum^{I_N}_{i_n=1} \mathcal{X}_{i_1,i_2,\dots,i_N}^2)^{1/2},
\end{equation}
where $\mathcal{X}_{i_1,i_2,\dots,i_N}$ means a tensor entry with the coordinate of $i_1,i_2,\dots,i_N$.
\end{definition}

\begin{definition} [Tucker decomposition \cite{review1}]
	The Tucker decomposition of an $N$th order tensor $\mathcal{X}\in \mathbb{R}^{I_{1} \times I_{2} \times \cdots \times I_{N}}$ is defined as:
	\begin{equation}
	\mathcal{X}\approx\mathcal{S}\times_1 U_1 \times_2 U_2 \cdots \times_n U_n	\cdots \times_N U_N,
	\end{equation}
	where tensor $\mathcal{S}\in\mathbb{R}^{r_1\times r_2\times \cdots r_N}$ is named as the core tensor, and matrix	$U_n\in\mathbb{R}^{I_n\times r_n}$, $n=1,2,\dots,N$, is called the factor matrix.
\end{definition}

\begin{definition} [Tensor mode-$n$ unfolding]
The mode-$n$ unfolding of an $N$th order tensor $\mathcal{X}\in \mathbb{R}^{I_{1} \times I_{2} \times \cdots \times I_{N}}$ produces a matrix $X\in \mathbb{R}^{I_n \times (\prod_{m=1, m\neq n}^N I_m)}$. In other words, a tensor entry $\mathcal{X}_{i_1,i_2,\dots,i_N}$ maps to a matrix entry ${X}_{i_n,j}$, where index follows
\begin{equation}
%\begin{aligned}
j=1+\sum_{k=1,k\neq n}^{N}(i_k-1)\prod_{m=1, m\neq n}^{k-1}I_m\text{.}
%J_k=\prod_{m=1, m\neq n}^{k-1}I_m.	
%\end{aligned}
\end{equation}
For brevity, the mode-$n$ unfolding of tensor $\mathcal{X}$ is represented as $\mathcal{X}_{(n)}$, and its inverse operation is denoted as $\operatorname{fold}_n(\cdot)$, namely, $\mathcal{X}=\operatorname{fold}_n(\mathcal{X}_{(n)})$.
\end{definition}

\begin{definition} [Tensor $n$-unfolding]
	The $n$-unfolding of an $N$th order tensor $\mathcal{X}\in \mathbb{R}^{I_{1} \times I_{2} \times \cdots \times I_{N}}$ leads to a  matrix $X\in \mathbb{R}^{(\prod_{m=1}^n I_m) \times (\prod_{m=n+1}^N I_m)}$. Furthermore, a tensor entry $\mathcal{X}_{i_1,i_2,\dots,i_N}$ maps to a matrix entry ${X}_{i,j}$, where index follows
	\begin{equation}
	\begin{aligned}
	&i=i_1+\sum_{k=2}^{n}(i_k-1)\prod_{m=1}^{k-1}I_m\text{,}\\ &j=i_{n+1}+\sum_{k=n+2}^{N}(i_k-1)\prod_{m=n+1}^{k-1}I_m{.}
	%J_k=\prod_{m=1}^{k-1}I_m.
	\end{aligned}
	\end{equation}
	The $n$-unfolding is a balanced tensor unfolding scheme.
\end{definition}

\begin{definition} [Kronecker product \cite{kron-p,kr-p}]
	The Kronecker product is denoted as $\otimes$. For two arbitrary matrices $X\in \mathbb{R}^{I_1\times I_2}$, $Y\in \mathbb{R}^{J_1\times J_2}$, we have
	\begin{equation}
		Z=X\otimes Y \in \mathbb{R}^{(I_1J_1)\times(I_2J_2)},
	\end{equation}
	where the matrix entry $Z_{k_1,k_2}=X_{i_1,i_2}Y_{j_1,j_2}$, $k_1=j_1+(i_1-1)J_1$, $k_2=j_2+(i_2-1)J_2$. Particularly, the Kronecker product admits the column-wise form, which is denoted as $\odot$. For matrices with same column size $X\in \mathbb{R}^{I\times K}$, $Y\in \mathbb{R}^{J\times K}$, we have
	\begin{equation}
	Z=X\odot Y \in \mathbb{R}^{(IJ)\times K}.
	\end{equation}
	The column-wise Kronecker product is also named as the Khatri-Rao product.
\end{definition}

\begin{definition} [Total variation \cite{tv-define}]
	The Total variation of matrix ${X}\in \mathbb{R}^{I_{1} \times I_{2}}$ consists of two formulations, namely, the isotropic TV and anisotropic TV. The isotropic TV of matrix $X$ is defined as:
	\begin{equation}
	\begin{aligned}
	\operatorname{TV}_i(X)=&(|\sum^{I_{1}-1}_{i=1}\sum^{I_{2}}_{j=1} X_{i+1,j}-X_{i,j}|^2\\
	&+|\sum^{I_{1}}_{i=1}\sum^{I_{2}-1}_{j=1} X_{i,j+1}-X_{i,j}|^2)^{1/2}.
	\end{aligned}
	\end{equation}
	The anisotropic TV of $X$ is defined as:
	\begin{equation}
	%\begin{aligned}
	\operatorname{TV}_a(X)=|\sum^{I_{1}-1}_{i=1}\sum^{I_{2}}_{j=1} X_{i+1,j}-X_{i,j}|+|\sum^{I_{1}}_{i=1}\sum^{I_{2}-1}_{j=1} X_{i,j+1}-X_{i,j}|,
	%\end{aligned}
	\end{equation}
The total variation can be generalized to the tensor form.	
\end{definition}

\begin{definition} [Tensor inner product]
	The inner product of two tensors $\mathcal{X}$, $\mathcal{Y} \in \mathbb{R}^{I_1\times I_2\times \cdots \times I_N}$ is defined as:
	\begin{equation}
		<\mathcal{X},\mathcal{Y}>=\sum_{i_1,i_2,\dots,i_N}\mathcal{X}_{i_1,i_2,\dots,i_N}\mathcal{Y}_{i_1,i_2,\dots,i_N}.
	\end{equation}
	The two tensors are with the same size to decide their inner product.
\end{definition}

\subsection{Related Works}
This subsection shows the related works on the low-rank prior and the smoothness prior of tensor completion. The low-rank tensor completion for incomplete tensor $\mathcal{X}\in \mathbb{R}^{I_{1} \times I_{2} \times \cdots \times I_{N}}$ is formulated as an optimization problem:
\begin{equation}
\label{lrtc-ori}
	\min_{\mathcal{X}}\ \operatorname{rank}(\mathcal{X}) \text{\ \ s.t.\ } P_\Omega(\mathcal{X})=P_\Omega(\mathcal{T}),
\end{equation}
where $\Omega$ suggests the observation set, which tells the cooridinates of given tensor entries. $P$ is the operator, mapping the tensor entries from original tensor $\mathcal{T}$ to incomplete tensor $\mathcal{X}$. 

Since the problem in equation (\ref{lrtc-ori}) is NP-hard, Zhao \textit{et al.} \cite{cp4} applied Bayesian inference to decide an appropriate tensor CP rank, but the estimated tensor CP rank may be inaccurate. Furthermore, Liu \textit{et al.}  \cite{tu3} used tensor mode-$n$ unfolding, transforming the low-rank tensor completion to the low-rank matrix completion: 
\begin{equation}
\label{lrtc-ha}
	\min_{\mathcal{X}_{(n)}}\ \sum_{n=1}^{N}\alpha_n\|\mathcal{X}_{(n)}\|_* \text{\ \ s.t.\ } P_\Omega(\mathcal{X})=P_\Omega(\mathcal{T}),
\end{equation}
where $\|\mathcal{X}_{(n)}\|_*$ decides the matrix nuclear norm of mode-$n$ unfolding of tensor $\mathcal{X}$, and weight $\alpha_n$, $n=1,2,\dots,N$, follows $\sum_{n=1}^{N}\alpha_n=1$. Considering the computational complexity of performing the matrix nuclear norm, Xu \textit{et al.} \cite{tu2} employed the parallel matrix factorization to reformulate the problem in equation (\ref{lrtc-ha}) as:
\begin{equation}
\label{lrtc-tmac}
\min_{\mathcal{X}_{(n)}, U_n, V_n}\ \sum_{n=1}^{N}\alpha_n\|\mathcal{X}_{(n)}-U_nV_n\|^2_F \text{\ \ s.t.\ } P_\Omega(\mathcal{X})=P_\Omega(\mathcal{T}),
\end{equation}
where $U_n\in \mathbb{R}^{I_n\times r_n}$, $V_n\in \mathbb{R}^{r_n\times (\prod^N_{m=1, m\neq n}I_m)}$, are factor matrices to $\mathcal{X}_{(n)}$, and $r_n$ is a prespecified low-rank prior. Both the methods efficiently conduct the tensor completion of incomplete tensors, but the tensor unfolding destroys the global structure. For remedying the structure damage, Zhang \textit{et al.} \cite{tubal2} employed the matrix singular value decomposition in Fourier domain to gain the low-rank prior. Although global structure is well preserved, it is only applicable to third-order tensors. Additionally, Bengua \textit{et al.} \cite{tt3} utilized $n$-unfolding to replace mode-$n$ unfolding, which better preserves the global structure. However, the $n$-unfolding produces unbalanced border matrices, since the border ranks are fixed.

The tensor low-rank priors focus on the global structure, while the smoothness priors are crucial for the local pattern of tensors. Chen \textit{et al.} \cite{tv3} simultaneously considered the low rank prior and the smoothness prior. Based on the Tucker decomposition, they proposed the model as:
\begin{equation}
\label{stdc}
\begin{aligned}
&\min_{\mathcal{X},\mathcal{S},U_n}\ \sum_{n=1}^{N}\alpha_n\|{U}_{n}\|_*\\
&+\beta \operatorname{tr}((U_1\otimes\cdots \otimes U_N)L(U_1\otimes\cdots \otimes U_N)^T)+\gamma \|\mathcal{S}\|^2_F\\
&\text{s.t.\ }\mathcal{X}=\mathcal{S}\times_1 U_1 \cdots \times_N U_N,\ P_\Omega(\mathcal{X})=P_\Omega(\mathcal{T}),
\end{aligned}
\end{equation}
where $\beta$, $\gamma$ are regularization parameters, $\operatorname{tr}(\cdot)$ is the trace operation, and matrix $L\in \mathbb{R}^{(\prod_{n=1}^N I_n)\times (\prod_{n=1}^N I_n)}$ is Laplacian matrix. The model in equation (\ref{stdc}) can be interpreted via two parts, i.e., the nuclear norm term for the low-rank prior, the Laplacian term and Frobenius norm term for the smoothness prior. In spite of this comprehensive using of both the low-rank prior and the smoothness prior, the computational complexity may reduce the performance. 

For further improvement, Ji \textit{et al.} \cite{mftv} fused the TV prior to the problem in equation (\ref{lrtc-tmac}) as follows:
\begin{equation}
\label{lrtc-mftv}
\begin{aligned}
&\min_{\mathcal{X}_{(n)}, U_n, V_n}\ \sum_{n=1}^{N}\alpha_n\|\mathcal{X}_{(n)}-U_nV_n\|^2_F+\mu \operatorname{TV}_i (V_3),\\
&\text{\ \ s.t.\ } P_\Omega(\mathcal{X})=P_\Omega(\mathcal{T}),
\end{aligned}
\end{equation}
where $\operatorname{TV}_i(\cdot)$ indicates the isotropic TV calculation, $\mu$ is the compromise parameter. Concretely, the matrix $U_n$ represents the library, and $V_n$ is treated as the unmixing of $\mathcal{X}_{(n)}$. Since they claimed the mode-$3$ of the tensor might contain the complete information and enjoy the useful structure, the TV prior is only applied to $V_3$. Besides, Jiang \textit{et al.} \cite{rank+smooth3} designed the framelet transform based dictionary to provide SC prior for $V_3$, which shows obvious progress in smoothness performance.  Nevertheless, the local pattern may not be fully exploited due to the partial utilization of the mode-$3$ unfolding. 

\section{Model}
\label{sec-model}
To sufficiently exploit both the global low-rank prior and the local smoothness prior, we first reformulate the tensor mode-$n$ unfolding, $n=1,2,\dots,N$, as follows:
\begin{equation}
	\mathcal{X}_{(n)}\approx U_n(U_N\odot\cdots \odot U_{n+1}\odot U_{n-1}\odot \cdots \odot U_{1})^T,
\end{equation}
where $U_n\in\mathbb{R}^{I_n\times r_n}$, indicates the matrix with a  unified low-rank prior $r_n$ of $n$th tensor mode, and $\odot$ implies the Khatri-Rao product. In addition, we denote $V_n=(U_N\odot\cdots \odot U_{n+1}\odot U_{n-1}\odot \cdots \odot U_{1})^T$, $V_n\in \mathbb{R}^{r_n\times (\prod^N_{m=1, m\neq n}I_m)}$ for notational brevity. 

\begin{figure}[t]
	\centering
	\begin{center}%[f]{0.25\linewidth}
		%\text{{\tiny NRSE}}\vspace{-0.03cm}	
		\includegraphics[width=8.6cm]{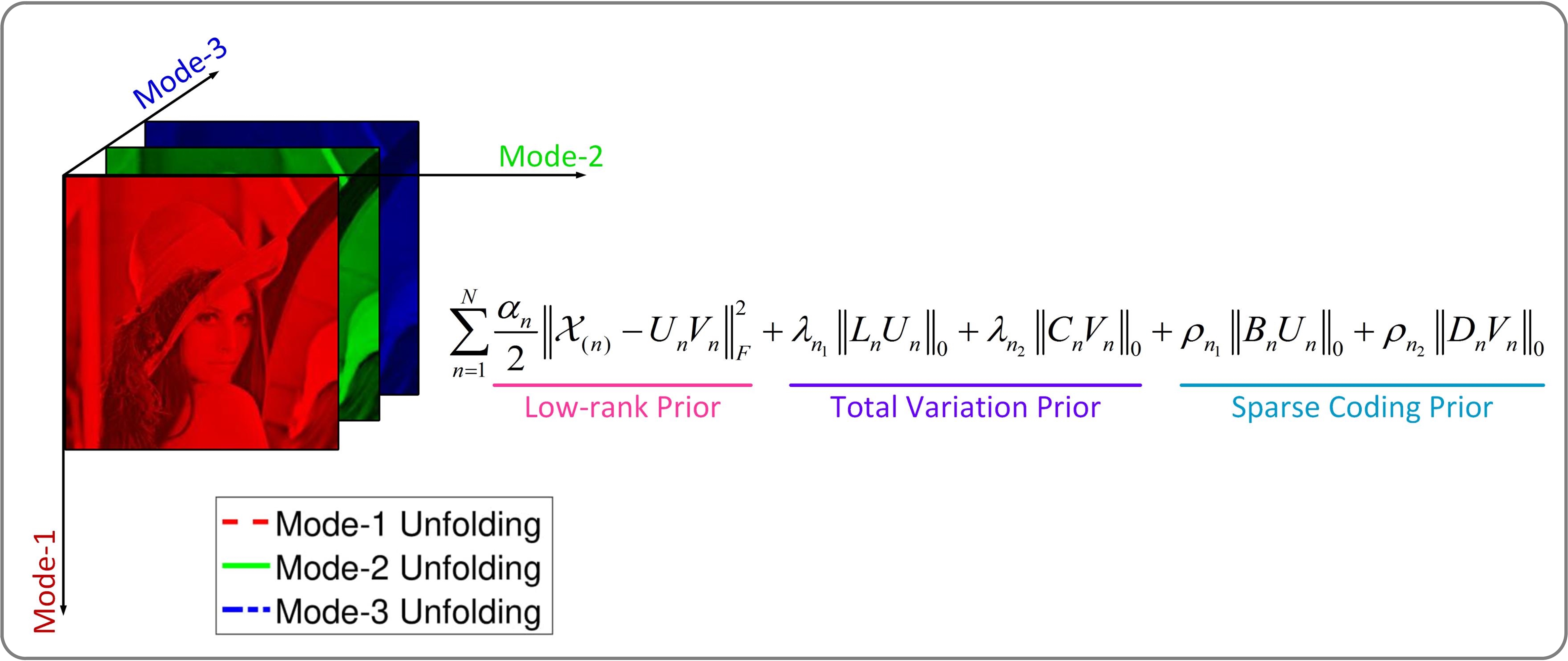}
		\\
		\text{{\footnotesize (a)}}%
	\end{center}%
	\begin{minipage}[f]{0.332\linewidth}
		\centering
		\includegraphics[height=3.5cm,width=2.74cm]{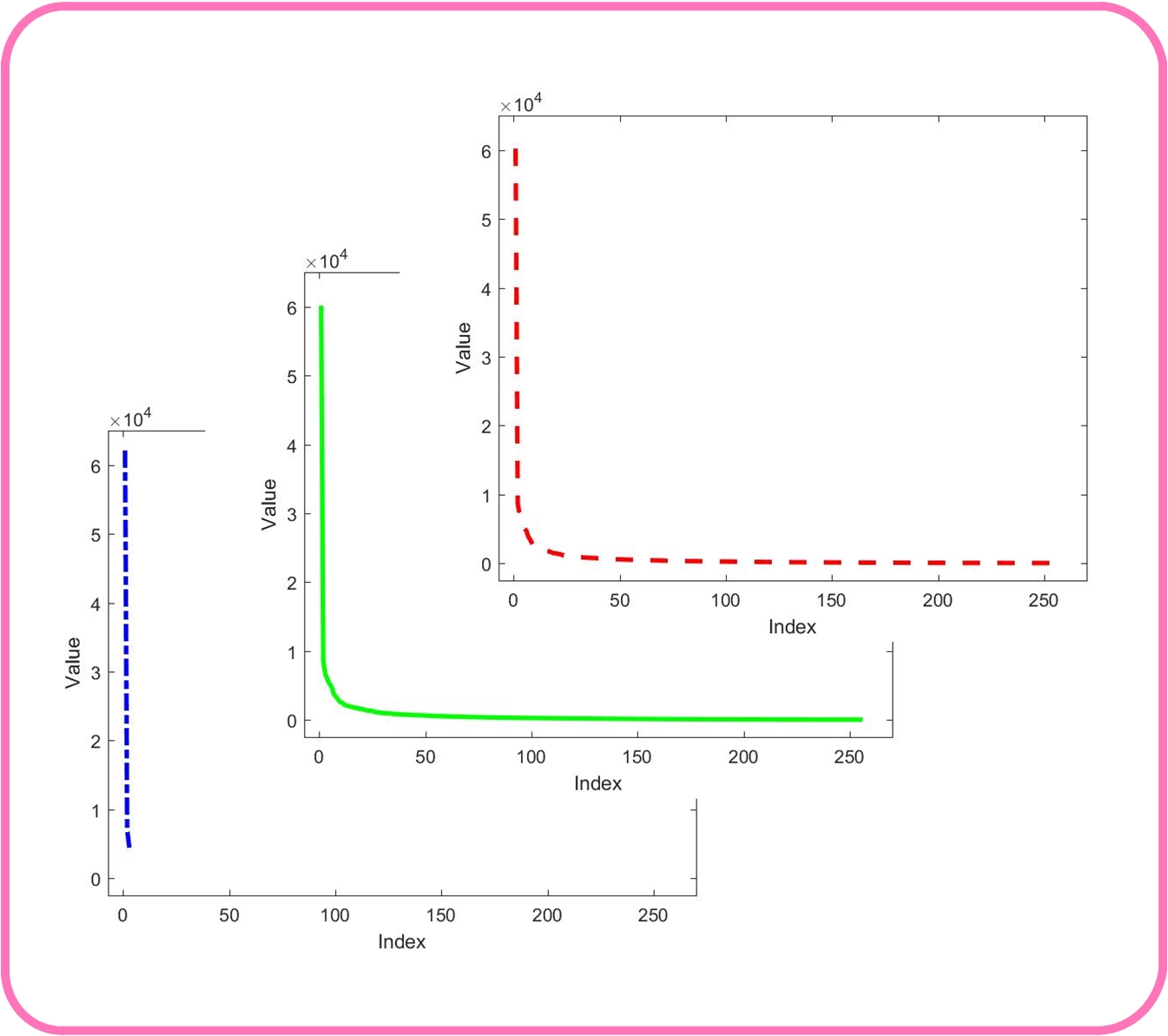}%
		\\
		\text{{\footnotesize (b)}}
	\end{minipage}%
	\begin{minipage}[f]{0.332\linewidth}
		\centering
		\includegraphics[height=3.5cm,width=2.74cm]{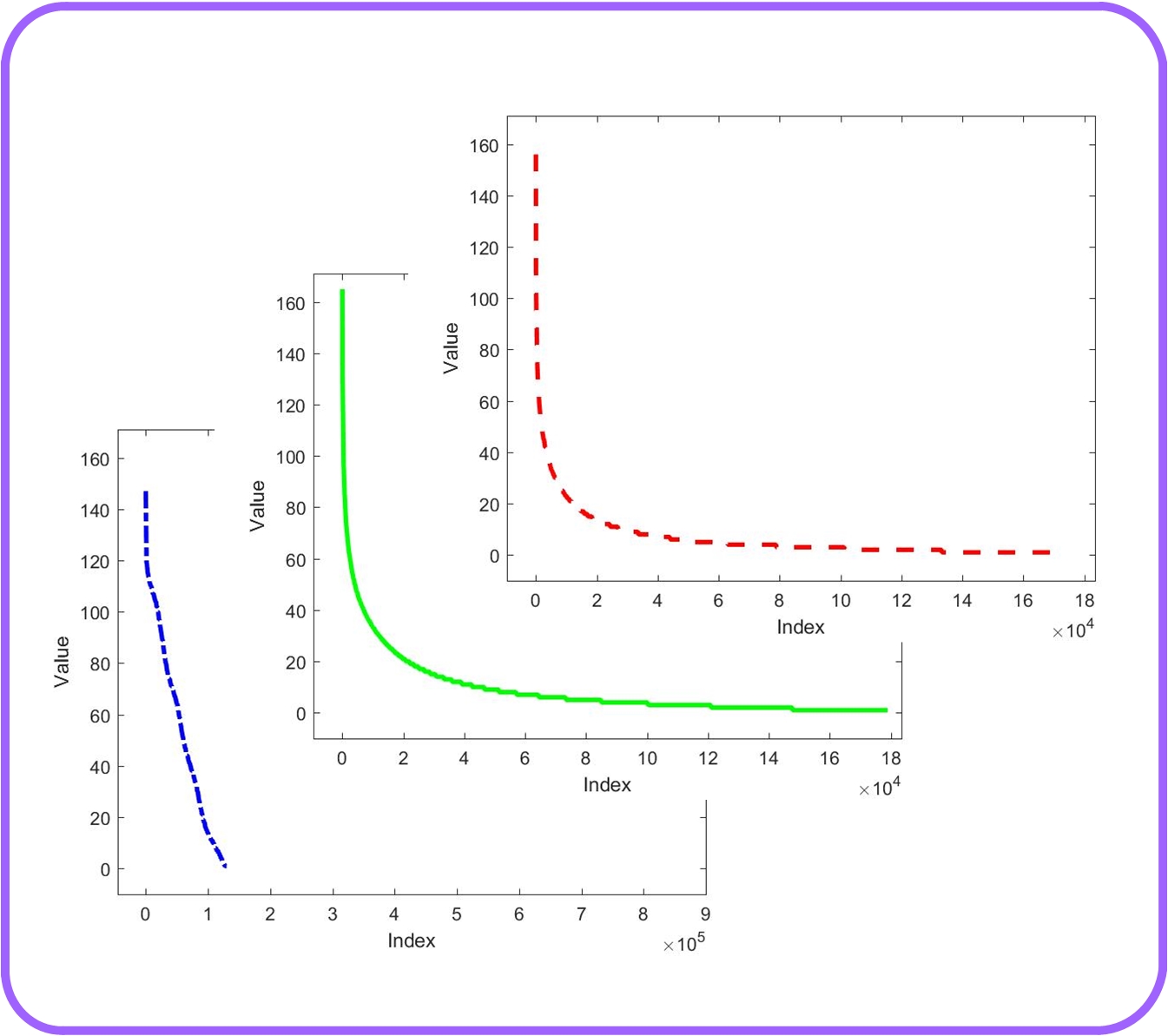}%
		\\
		\text{{\footnotesize (c)}}
	\end{minipage}%
	\begin{minipage}[f]{0.332\linewidth}
		\centering
		\includegraphics[height=3.5cm,width=2.74cm]{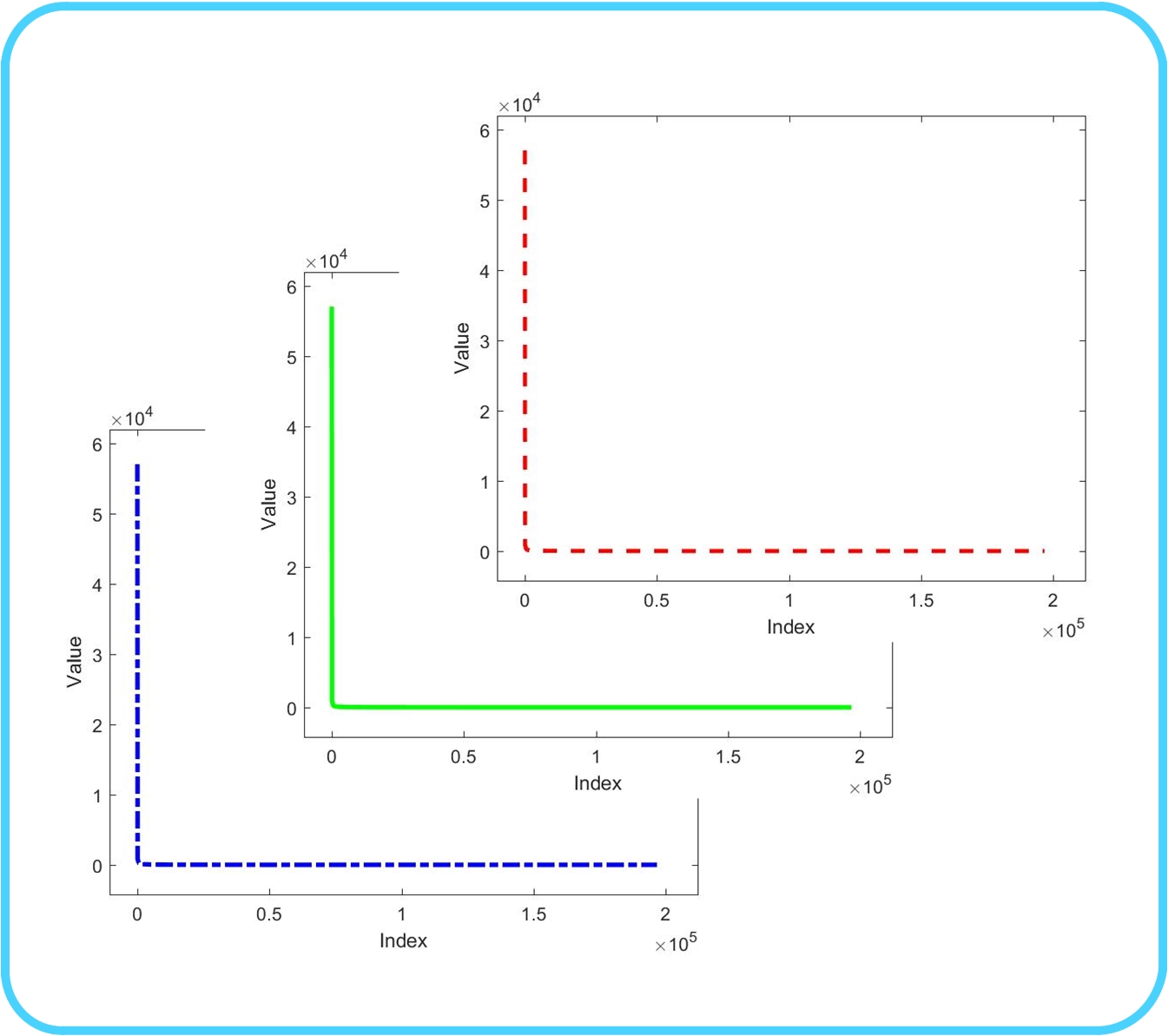}%
		\\
		\text{{\footnotesize (d)}}
	\end{minipage}%
	\caption{Graphical explaination of the proposed model. (a) The hierarchical prior regularized matrix factorization model for a color image $\mathcal{X}$ of size ${256\times 256\times 3}$. (b) First hierarchy for low-rank prior: Singular values of mode-$n$ unfoldings. (c) Second hierarchy for TV prior: Margin values of mode-$n$ unfoldings. (d) Third hierarchy for SC prior: Absolute coefficient values of mode-$n$ unfoldings by DCT.}
	\label{graph-illus}
\end{figure}

We consider each tensor mode-$n$ unfolding can be hierarchically characterized via the low-rank prior, TV prior, and SC prior. Accordingly, the following optimization model is established:
	\begin{equation}
	\label{proposed-model}
	\begin{aligned}
	\min_{\mathcal{X}_{(n)}, U_n, V_n}\ &\sum_{n=1}^{N}\frac{\alpha_n}{2}\|\mathcal{X}_{(n)}-U_nV_n\|^2_F\\
	&+\lambda_{n_1}\|L_nU_n\|_0+\lambda_{n_2}\|C_nV_n\|_0\\
	&+\rho_{n_1}\|B_nU_n\|_0+\rho_{n_2}\|D_nV_n\|_0,\\
	&\text{s.t.\ } P_\Omega(\mathcal{X})=P_\Omega(\mathcal{T}),
	\end{aligned}
	\end{equation}
where $L_n\in \mathbb{R}^{(I_n-1)\times I_n}$ and $C_n\in \mathbb{R}^{(r_n-1)\times r_n}$ represent the TV regularization matrices, e.g., entries $L_{i,i}=1$, $L_{i,i+1}=-1$, $i=1,2,\dots,I_n-1$, for matrix $U_n$, while $C_n$ is for matirx $V_n$. Matrices $B_n\in\mathbb{R}^{I_n\times I_n}$ and $D_n\in\mathbb{R}^{r_n\times r_n}$ denote the DCT dictionaries for $U_n$, $V_n$ respectively. The weight $\alpha_n$ observes $\sum_{n=1}^{N}\alpha_n=1$. Compromising parameters are shown as $\lambda_{n_1}$, $\lambda_{n_2}$, $\rho_{n_1}$, $\rho_{n_2}$, and $\|\cdot\|_0$ suggests the matrix $l_0$ norm. For example, Fig. \ref{graph-illus} explains the proposed hierarchical prior regularized matrix factorization model for a color image.

Unfortunately, the optimization problem in equation (\ref{proposed-model}) is NP-hard owing to the matrix $l_0$ norm. Hence, we use the convex surrogate to further reformulate the problem into:
	\begin{equation}
\label{proposed-model2}
\begin{aligned}
\min_{\mathcal{X}_{(n)}, U_n, V_n}\ &\sum_{n=1}^{N}\frac{\alpha_n}{2}\|\mathcal{X}_{(n)}-U_nV_n\|^2_F\\
&+\lambda_{n_1}\|L_nU_n\|_1+\lambda_{n_2}\|C_nV_n\|_1\\
&+\rho_{n_1}\|B_nU_n\|_1+\rho_{n_2}\|D_nV_n\|_1,\\
&\text{s.t.\ } P_\Omega(\mathcal{X})=P_\Omega(\mathcal{T}),
\end{aligned}
\end{equation}
where $\|\cdot\|_1$ determines the matrix $l_1$ norm. Introducing auxiliary variables $G_n$, $H_n$, $R_n$, $M_n$, $n=1,2,\dots,N$, the problem in equation (\ref{proposed-model2}) is equivalent to:
\begin{equation}
\label{proposed-model3}
\begin{aligned}
\min_{\substack{\mathcal{X}_{(n)}, U_n, V_n,\\ G_n, H_n, R_n, M_n}}&\sum_{n=1}^{N}\frac{\alpha_n}{2}\|\mathcal{X}_{(n)}-U_nV_n\|^2_F\\
&+\lambda_{n_1}\|G_n\|_1+\lambda_{n_2}\|H_n\|_1\\
&+\rho_{n_1}\|R_n\|_1+\rho_{n_2}\|M_n\|_1,\\
&\text{s.t.\ } G_n=L_nU_n,H_n=C_nV_n,\\
&R_n=B_nU_n,M_n=D_nV_n,\\ 
&P_\Omega(\mathcal{X})=P_\Omega(\mathcal{T}).
\end{aligned}
\end{equation}

The model in equation (\ref{proposed-model3}) can be parallelly solved for each mode. Therefore, we establish the corresponding augmented Lagrangian function $f_{L}$ to the $n$th mode of the tensor as:
\begin{equation}
\label{proposed-model4}
\begin{aligned}
&f_L({\mathcal{X}_{(n)}, U_n, V_n, G_n, H_n, R_n, M_n,\Lambda_n, \Pi_n, \Phi_n,  \Gamma_n})\\
&=\frac{\alpha_n}{2}\|\mathcal{X}_{(n)}-U_nV_n\|^2_F\\
&+\lambda_{n_1}\|G_n\|_1+<\Lambda_n,L_nU_n-G_n>+\frac{\beta_{n_1}}{2}\|L_nU_n-G_n\|_F^2\\
&+\lambda_{n_2}\|H_n\|_1+<\Pi_n,C_nV_n-H_n>+\frac{\beta_{n_2}}{2}\|C_nV_n-H_n\|_F^2\\
&+\rho_{n_1}\|R_n\|_1+<\Phi_n,B_nU_n-R_n>+\frac{\omega_{n_1}}{2}\|B_nU_n-R_n\|_F^2\\
&+\rho_{n_2}\|M_n\|_1+<\Gamma_n,D_nV_n-M_n>+\frac{\omega_{n_2}}{2}\|D_nV_n-M_n\|_F^2\\
&\text{s.t.\ } P_\Omega(\mathcal{X})=P_\Omega(\mathcal{T}),
\end{aligned}
\end{equation}
where $\Lambda_n$, $\Pi_n$, $\Phi_n$ $\Gamma_n$ are dual variables to the equality constraints, and $\beta_{n_1}$, $\beta_{n_2}$, $\omega_{n_1}$, $\omega_{n_2}$ are augmented Lagrangian parameters for the $n$th mode of the tensor. The formulation in equation (\ref{proposed-model4}) are subsequently rewrote as:
\begin{equation}
\label{proposed-model5}
\begin{aligned}
&f_L({\mathcal{X}_{(n)}, U_n, V_n, G_n, H_n, R_n, M_n,\Lambda_n, \Pi_n, \Phi_n,  \Gamma_n})\\
&=\frac{\alpha_n}{2}\|\mathcal{X}_{(n)}-U_nV_n\|^2_F\\
&+\lambda_{n_1}\|G_n\|_1+\frac{\beta_{n_1}}{2}\|L_nU_n-G_n+\frac{\Lambda_{n}}{\beta_{n_1}}\|_F^2\\
&+\lambda_{n_2}\|H_n\|_1+\frac{\beta_{n_2}}{2}\|C_nV_n-H_n+\frac{\Pi_n}{\beta_{n_2}}\|_F^2\\
&+\rho_{n_1}\|R_n\|_1+\frac{\omega_{n_1}}{2}\|B_nU_n-R_n+\frac{\Phi_n}{\omega_{n_1}}\|_F^2\\
&+\rho_{n_2}\|M_n\|_1+\frac{\omega_{n_2}}{2}\|D_nV_n-M_n+\frac{\Gamma_n}{\omega_{n_2}}\|_F^2 
+\Delta\\
&\text{s.t.\ } P_\Omega(\mathcal{X})=P_\Omega(\mathcal{T}),
\end{aligned}
\end{equation}
where $\Delta$ is a constant. 

\section{Algorithm}
\label{sec-algo}
This section demonstrates the algorithm for the proposed model. The convergence and complexity of the algorithm are also investigated. 

\subsection{Solution for Algorithm}
We address the problem in equation (\ref{proposed-model5}) under the framework of ADMM. For the $k$th iteration, we have the following main subproblems:

\begin{enumerate}[fullwidth,itemindent=1em]
\item $U_n$ \textbf{subproblem}: We update $U_n$ at $k$th iteration as 
\begin{equation}
\label{subp1}
\begin{aligned}
U^{(k+1)}_n
&=\arg \min_{U_n} \frac{\alpha_n}{2}\|\mathcal{X}^{(k)}_{(n)}-U_nV^{(k)}_n\|^2_F\\
&+\frac{\beta^{(k)}_{n_1}}{2}\|L_nU_n-G^{(k)}_n+\frac{\Lambda^{(k)}_{n}}{\beta^{(k)}_{n_1}}\|_F^2\\
&+\frac{\omega^{(k)}_{n_1}}{2}\|B_nU_n-R^{(k)}_n+\frac{\Phi^{(k)}_n}{\omega^{(k)}_{n_1}}\|_F^2.
\end{aligned}
\end{equation}
This formulation admits a closed form solution. We then derive that:
\begin{equation}
\begin{aligned}
	&\alpha_n U_n V^{(k)}_nV_n^{{(k)},T}+(\beta^{(k)}_{n_1} L_n^TL_n+\omega^{(k)}_{n_1} B_n^TB_n) U_n\\ 
	&=\alpha_n X_{(n)}V_n^{{(k)},T}+\beta^{(k)}_{n_1}L_n^TG^{(k)}_n-L_n^T\Lambda^{(k)}_n\\
	&+\omega^{(k)}_{n_1} B_n^TR_n^{(k)}-B_n^T\Phi^{(k)}_n,
\end{aligned}
\end{equation}
where $V_n^{{(k)},T}$ represents the transpose of $V_n^{{(k)}}$. By the vectorization operation $\operatorname{vec}(\cdot)$, we have:
\begin{equation}
\begin{aligned}
&(\alpha_n V^{(k)}_nV_n^{{(k)},T} \otimes I_{a}\\
&+ I_{b} \otimes (\beta^{(k)}_{n_1} L_n^TL_n+\omega^{(k)}_{n_1} B_n^TB_n))\operatorname{vec}(U_n)\\ 
&=\operatorname{vec}(\alpha_n X_{(n)}V_n^{{(k)},T}+\beta^{(k)}_{n_1}L_n^TG^{(k)}_n-L_n^T\Lambda^{(k)}_n\\
&+\omega^{(k)}_{n_1} B_n^TR_n^{(k)}-B_n^T\Phi^{(k)}_n
),
\end{aligned}
\end{equation}
where $I_{a}\in \mathbb{R}^{I_n\times I_n}$, $I_{b}\in \mathbb{R}^{r_n\times r_n}$ represent the identity matrices, $\otimes$ is the Kronecker product, and it provides:
\begin{equation}
\label{admm-u}
\begin{aligned}
& U^{(k+1)}_n = \operatorname{ivec}(
(\alpha_n V^{(k)}_nV_n^{{(k)},T} \otimes I_{a}\\
&+I_{b} \otimes (\beta^{(k)}_{n_1} L_n^TL_n+\omega^{(k)}_{n_1} B_n^TB_n))^{\dagger}\\
& \operatorname{vec}(\alpha_n X_{(n)}V_n^{{(k)},T}+\beta^{(k)}_{n_1}L_n^TG^{(k)}_n-L_n^T\Lambda^{(k)}_n\\
&+\omega^{(k)}_{n_1} B_n^TR_n^{(k)}-B_n^T\Phi^{(k)}_n)),
\end{aligned}
\end{equation}
where $\operatorname{ivec}(\cdot)$ is the inverse operation of $\operatorname{vec}(\cdot)$, and $\dagger$ denotes the Moore-Penrose pseudo-inverse.
\item $V_n$ \textbf{subproblem}: At $k$th iteration, we update matrix $V_n$ via
\begin{equation}
\label{subp2}
\begin{aligned}
&V^{(k+1)}_n
=\arg \min_{V_n} \frac{\alpha_n}{2}\|\mathcal{X}^{(k)}_{(n)}-U^{(k+1)}_nV_n\|^2_F\\
&+\frac{\beta^{(k)}_{n_2}}{2}\|C_nV_n-H^{(k)}_n+\frac{\Pi^{(k)}_n}{\beta^{(k)}_{n_2}}\|_F^2\\
&+\frac{\omega^{(k)}_{n_2}}{2}\|D_nV_n-M^{(k)}_n+\frac{\Gamma^{(k)}_n}{\omega^{(k)}_{n_2}}\|_F^2.
\end{aligned}
\end{equation}
Owing to the convexity this subproblem, we consequently present that:
\begin{equation}
\label{admm-v}
\begin{aligned}
&V^{(k+1)}_n
=(\alpha_nU_n^{(k+1),T}U^{(k+1)}_n+\beta^{(k)}_{n_2}C_n^TC_n+\omega^{(k)}_{n_2}I_b)^{\dagger}\\
&(\alpha_nU_n^{(k+1),T}\mathcal{X}^{(k)}_{(n)}+\beta^{(k)}_{n_2}C_n^TH^{(k)}_n-C_n^T\Pi^{(k)}_n\\
&+\omega^{(k)}_{n_2}D_n^TM^{(k)}_n-D_n^T\Gamma^{(k)}_n),
\end{aligned}
\end{equation}
where $I_{b}\in \mathbb{R}^{r_n\times r_n}$ is the identity matrix.
\item $G_n$ \textbf{Subproblem}: Matrix $G_n$ is updated at $k$th iteration by solving
\begin{equation}
\label{subp3}
\begin{aligned}
G^{(k+1)}_n
&=\arg \min_{G_n} \lambda_{n_1}\|G_n\|_1\\
&+\frac{\beta^{(k)}_{n_1}}{2}\|L_nU^{(k+1)}_n-G_n+\frac{\Lambda^{(k)}_n}{\beta^{(k)}_{n_1}}\|_F^2.
\end{aligned}
\end{equation}
The matrix $l_1$ norm can be solved via the soft value thresholding operation $\tau_\epsilon(\cdot)$ \cite{svt}. Subsequently, we have:
\begin{equation}
\label{admm-g}
	G^{(k+1)}_n=\tau_{\frac{\lambda_{n_1}}{\beta^{(k)}_{n_1}}}(L_nU^{(k+1)}_n+\frac{\Lambda^{(k)}_n}{\beta^{(k)}_{n_1}}).
\end{equation}
\item $H_n$ \textbf{subproblem}: Matrix $H_n$ is updated at $k$th iteration by solving
\begin{equation}
\label{subp4}
\begin{aligned}
H^{(k+1)}_n
&=\arg \min_{H_n} \lambda_{n_2}\|H_n\|_1\\
&+\frac{\beta^{(k)}_{n_2}}{2}\|C_nV^{(k+1)}_n-H_n+\frac{\Pi^{(k)}_n}{\beta^{(k)}_{n_2}}\|_F^2.
\end{aligned}
\end{equation}
Similarly, we have:
\begin{equation}
\label{admm-h}
H^{(k+1)}_n=\tau_{\frac{\lambda_{n_2}}{\beta^{(k)}_{n_2}}}(C_nV^{(k+1)}_n+\frac{\Pi^{(k)}_n}{\beta^{(k)}_{n_2}}).
\end{equation}
\item $R_n$ \textbf{subproblem}: Matrix $R_n$ is updated at $k$th iteration by solving
\begin{equation}
\label{subp5}
\begin{aligned}
R^{(k+1)}_n
&=\arg \min_{R_n} \rho_{n_1}\|R_n\|_1\\
&+\frac{\omega^{(k)}_{n_1}}{2}\|B_nU^{(k+1)}_n-R_n+\frac{\Phi^{(k)}_n}{\omega^{(k)}_{n_1}}\|_F^2.
\end{aligned}
\end{equation}
Consequently, we have:
\begin{equation}
\label{admm-r}
R^{(k+1)}_n=\tau_{\frac{\rho_{n_1}}{\omega^{(k)}_{n_1}}}(B_nU^{(k+1)}_n+\frac{\Phi^{(k)}_n}{\omega^{(k)}_{n_1}}).
\end{equation}
\item $M_n$ \textbf{subproblem}: Matrix $M_n$ is updated at $k$th iteration by solving
\begin{equation}
\label{subp6}
\begin{aligned}
M^{(k+1)}_n
&=\arg \min_{M_n} \rho_{n_2}\|M_n\|_1\\
&+\frac{\omega^{(k)}_{n_2}}{2}\|D_nV^{(k+1)}_n-M_n+\frac{\Gamma^{(k)}_n}{\omega^{(k)}_{n_2}}\|_F^2.
\end{aligned}
\end{equation}
The solution to this problem is:
\begin{equation}
\label{admm-m}
M^{(k+1)}_n=\tau_{\frac{\rho_{n_2}}{\omega^{(k)}_{n_2}}}(D_nV^{(k+1)}_n+\frac{\Gamma^{(k)}_n}{\omega^{(k)}_{n_2}}).
\end{equation}
\item $\mathcal{X}_{(n)}$ \textbf{subproblem}: We update the tensor mode-$n$ unfolding at $k$th iteration as
	\begin{equation}
	\label{subp7}
	\mathcal{X}^{(k+1)}_{(n)}
	=\arg \min_{\mathcal{X}_{(n)}} \frac{\alpha_n}{2}\|\mathcal{X}_{(n)}-U^{(k+1)}_nV^{(k+1)}_n\|^2_F.
	\end{equation}
This subproblem has a closed form solution due to the convexity of the formulation:
\begin{equation}
\label{admm-x}
\mathcal{X}^{(k+1)}_{(n)}=U^{(k+1)}_nV^{(k+1)}_n.
\end{equation}
\item Multiplier \textbf{subproblem}: Dual variables $\Lambda_n$, $\Pi_n$, $\Phi_n$, $\Gamma_n$, are updated following the ADMM procedure, and we have
\begin{equation}
\label{admm-dual}
	\left\{
	\begin{aligned}
	&\Lambda^{(k+1)}_n=\Lambda^{(k)}_n+\beta_{n_1}(L_nU^{(k+1)}_n-G^{(k+1)}_n),\\
	&\Pi^{(k+1)}_n=\Pi^{(k)}_n+\beta_{n_2}(C_nV^{(k+1)}_n-H^{(k+1)}_n),\\
	&\Phi^{(k+1)}_n=\Phi^{(k)}_n+\omega_{n_1}(B_nU^{(k+1)}_n-R^{(k+1)}_n),\\
	&\Gamma^{(k+1)}_n=\Gamma^{(k)}_n+\omega_{n_2}(D_nV^{(k+1)}_n-M^{(k+1)}_n).
	\end{aligned}
	\right.
\end{equation}
\end{enumerate}
Moreover, we introduce a constant $\mu$ to accelerate the convergence, e.g., $\beta^{(k+1)}_{n_1}=\mu\beta^{(k)}_{n_1}$.

Finally, to employ the constraint provided by the observation set, we update the recovered tensor $\mathcal{X}$ using all mode-$n$ unfoldings as follows:
\begin{equation}
\label{admm-tensor}
	\mathcal{X}^{(k+1)}=P_\Omega(\mathcal{T})+P_{\Omega^{\perp}}(\sum_{n=1}^{N} \alpha_n \operatorname{fold}_n(\mathcal{X}^{(k+1)}_{(n)})),
\end{equation}
where $\Omega^{\perp}$ means the observation set of missing components. Furthermore, the convergence criterion of the proposed algorithm is to calculate the absolute relative error of two successive iterations as:
%\begin{equation}
%\left|\frac{\left\|\mathcal{X}^{(k+1)}\|_{F}-\|\mathcal{X}^{(k)}\right\|_{F}}{\left\|\mathcal{X}^{(k)}\right\|_{F}}\right| < \xi,	
%\end{equation}
\begin{equation}
\frac{\left|(\|\mathcal{X}^{(k+1)}\|_{F}-\|\mathcal{X}^{(k)}\|_{F})\right|}{\|\mathcal{X}^{(k)}\|_{F}} < \xi,	
\end{equation}
where $\xi$ is a small positive constant. Our method is summarized in Algorithm \ref{proposed-algorithm}.

\IncMargin{1em}
\begin{algorithm}
	\label{proposed-algorithm}
	\caption{Hierarchical Prior Regularized Matrix Factorization Algorithm (HPMF)}%算法名字
	\LinesNumbered %要求显示行号	
	\KwIn{Original $N$th-order tensor $\mathcal{T}$ with observation set $\Omega$. Parameters $\alpha_n$, $\lambda_{n_1}$, $\lambda_{n_2}$, $\rho_{n_1}$, $\rho_{n_2}$, $\beta^{(0)}_{n_1}$, $\beta^{(0)}_{n_2}$, $\omega^{(0)}_{n_1}$, $\omega^{(0)}_{n_2}$, constant $\mu$, maximum iterations $K$, TV regularization matrices $L_n$, $C_n$, DCT matrices $B_n$, $D_n$, $n=1,2,\dots,N$.}%输入参数	
	\textbf{Initialization:}  
	$P_\Omega(\mathcal{X}^{(0)})=P_\Omega(\mathcal{T})$, $P_{\Omega^{\perp
	}}(\mathcal{X}^{(0)})=0$, $U^{(0)}_n$, $V^{(0)}_n$, $k=0$\; 
		\While{not converged or $k<K$}
		{
			\For{$n=1,2,\dots,N$}
			{
			Update ${U}_n^{(k+1)}$ by equation (\ref{admm-u})\;	
			
			Update ${V}_n^{(k+1)}$ by equation (\ref{admm-v})\;
			
			Update ${G}_n^{(k+1)}$ by equation (\ref{admm-g})\;
			
			Update ${H}_n^{(k+1)}$ by equation (\ref{admm-h})\;
			
			Update ${R}_n^{(k+1)}$ by equation (\ref{admm-r})\;
			
			Update ${M}_n^{(k+1)}$ by equation (\ref{admm-m})\;
			
			Update $\mathcal{X}_{(n)}^{(k+1)}$ by equation (\ref{admm-x})\;
			
			Update $\Lambda^{(k+1)}_n$, $\Pi^{(k+1)}_n$, $\Phi^{(k+1)}_n$, $\Gamma^{(k+1)}_n$ by equation (\ref{admm-dual})\;
		}

			Update $\mathcal{X}^{(k+1)}$ by equation (\ref{admm-tensor})\;
			
			$\beta^{(k+1)}_{n_1}=\mu\beta^{(k)}_{n_1}$, 
			$\beta^{(k+1)}_{n_2}=\mu\beta^{(k)}_{n_2}$\;
			$\omega^{(k+1)}_{n_1}=\mu\omega^{(k)}_{n_1}$, 
			$\omega^{(k+1)}_{n_2}=\mu\omega^{(k)}_{n_2}$\;
			
			$k=k+1$.%\;
		}
		
	\KwOut{Recovered tensor $\mathcal{X}^{(k+1)}$.}% 输出
\end{algorithm}

\subsection{Complexity and Convergence}
We first discuss the overall complexity at each iteration of Algorithm \ref{proposed-algorithm}. For simplicity, we assume the mode lengths $I_n$ and the corresponding mode ranks $r_n$, $n=1,2,\dots,N$, ($N\geq3$), of an $N$th-order tensor are all equal to $I$ and $r$ respectively, $I\geq r$. Therefore, equation (\ref{admm-u}) shows the complexity $\mathcal{O}(r^3I^3+rI^N)$. Equation (\ref{admm-v}) shows the complexity $\mathcal{O}(rI^N)$. Equation (\ref{admm-g}) shows the complexity $\mathcal{O}(rI^2)$.  Equation (\ref{admm-h}) shows the complexity $\mathcal{O}(r^2I^{N-1})$. Equation (\ref{admm-r}) shows the complexity $\mathcal{O}(rI^2)$. Equation (\ref{admm-m}) shows the complexity $\mathcal{O}(r^2I^{N-1})$. Equation (\ref{admm-x}) shows the complexity $\mathcal{O}(rI^{N})$.  Equation (\ref{admm-dual}) shows the complexity $\mathcal{O}(r^2I^{N-1})$. Finally, equation (\ref{admm-tensor}) demonstrates the complexity $\mathcal{O}(I^{N})$. To sum up, the overall complexity of the proposed algorithm is $\mathcal{O}(r^3I^3+rI^N)$ at each iteration. Our algorithm can be parallelly solved along different tensor modes, which further improves the computational efficiency.

For convergence analysis, since the proposed model in equation (\ref{proposed-model5}) consists of separable convex subproblems, i.e., the squared $F$-norm and the $l_1$ norm subproblems, the convergence of Algorithm \ref{proposed-algorithm} is theoretically guaranteed under the ADMM framework \cite{admm-conver}.

\section{Experiment}
\label{sec-experi}
This section conducts experiments on color images. We compare the proposed algorithm with several state-of-the-art methods, including FBCP \cite{cp4}, HaLRTC \cite{tu3}, TNN \cite{tubal2}, SiLRTC-TT \cite{tt3}, STDC\cite{tv3}, and MF-TV \cite{mftv}. All experiments are implemented by the Intel i5-8500 CPU at 3.0 GHz, 32 GB RAM machine, under the Matlab (R2020a) environment.

We use sampling ratio (SR) to measure the incomplete tensors, and SR is defined as follows:
\begin{equation}
\text{SR}=\frac{|\Omega|}{\prod_{n=1}^{N} I_{n}},
\end{equation}
where $|\Omega|$ denotes the number of uniformly sampled entries at random from an original tensor $\mathcal{T}\in\mathbb{R}^{I_1\times I_2\times \cdots\times I_N}$. 

For numerical comparison, given the recovered tensor $\mathcal{X}$, we decide peak signal-to-noise ratio (PSNR) as:
\begin{equation}
\text{PSNR}=10 \log _{10} \big( \frac{\mathcal{T}_{\max}^{2}}{\|\mathcal{X}-\mathcal{T}\|_{F}^{2}}\cdot \prod_{n=1}^{N} I_{n}\big),
\end{equation}
where $\mathcal{T}_{\max}$ denotes the largest value of the original tensor $\mathcal{T}$. Relative standard error (RSE) is computed as:
\begin{equation}
	\text{RSE}=\frac{\|\mathcal{X}-\mathcal{T}\|_{F}}{\|\mathcal{T}\|_{F}}.
\end{equation}
For visual comparison, structure similarity (SSIM) assesses the similarity between $\mathcal{X}$ and $\mathcal{T}$ through luminance, contrast, and structure \cite{ssim}.

\subsection{Installation Settings}
For the proposed algorithm, the mode-$n$ unfolding rank $r_n$, $n=1,2,\dots,N$, is estimated as in the work \cite{tt3}, by keeping the singular values which satisfy:
\begin{equation}
	\frac{\sigma_i}{\sigma_1}>\delta,
\end{equation}
where $i=1,2,\dots,r_n$, and the singular value $\sigma_i$ is sorted in descending order. We emperically obtain the threshold $\delta$ from the range of $[0.01,0.37]$ for the best algorithm performance. Additionally, the convergence criterion $\xi$ is $10^{-5}$, and the maximum iterations $K$ is $500$. We decide the acceleration constant $\mu$ from the range of $[1,1.05]$ to obtain the best algorithm performance. The compared state-of-the-art algorithms utilize their default settings.

\subsection{Color Image with Random Sampling}
This part uses the typical third-order tensors, i.e., the color images, to test the recovering performance of algorithms. We consider the scenario that the experiments under different SRs. In this scenario, the compromising parameters of the proposed algorithm are set as $\alpha_n=1/3$, $\lambda_{n_1}=100$, $\lambda_{n_2}=100$, $\rho_{n_1}=0.1$, $\rho_{n_2}=100$. In addition, the augmented Lagrangian parameters are determined as $\beta^{(0)}_{n_1}=1$, $\beta^{(0)}_{n_2}=100$, $\omega^{(0)}_{n_1}=0.001$, $\omega^{(0)}_{n_2}=1000$ , where $n=1,2,3$.

\begin{figure}[htbp]
	\centering
	\begin{minipage}[f]{0.25\linewidth}
		\centering
		\includegraphics[width=2.15cm]{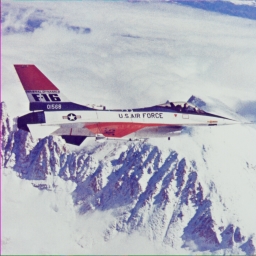}%
		\\
		\text{{\footnotesize (a)}}
	\end{minipage}%
	\begin{minipage}[f]{0.25\linewidth}
		\centering
		\includegraphics[width=2.15cm]{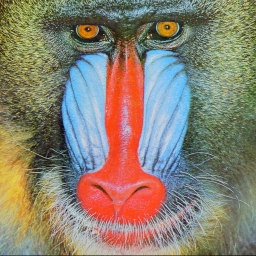}%
		\\
		\text{{\footnotesize (b)}}
	\end{minipage}%
	\begin{minipage}[f]{0.25\linewidth}
		\centering
		\includegraphics[width=2.15cm]{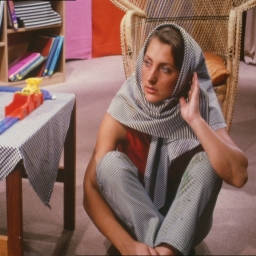}%
		\\
		\text{{\footnotesize (c)}}
	\end{minipage}%
	\begin{minipage}[f]{0.25\linewidth}
		\centering
		\includegraphics[width=2.15cm]{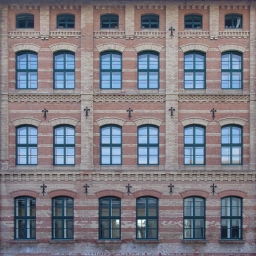}%
		\\
		\text{{\footnotesize (d)}}
	\end{minipage}%
\\%------
	\begin{minipage}[f]{0.25\linewidth}
		\centering
		\includegraphics[width=2.15cm]{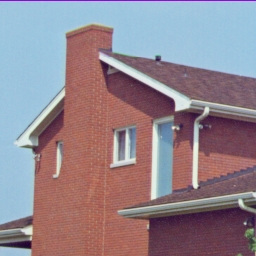}%
		\\
		\text{{\footnotesize (e)}}
	\end{minipage}%
	\begin{minipage}[f]{0.25\linewidth}
		\centering
		\includegraphics[width=2.15cm]{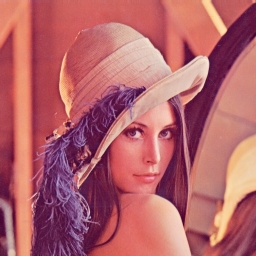}%
		\\
		\text{{\footnotesize (f)}}
	\end{minipage}%
	\begin{minipage}[f]{0.25\linewidth}
	\centering
	\includegraphics[width=2.15cm]{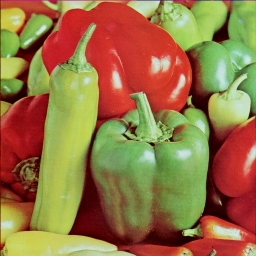}%
	\\
	\text{{\footnotesize (g)}}
\end{minipage}%
	\begin{minipage}[f]{0.25\linewidth}
	\centering
	\includegraphics[width=2.15cm]{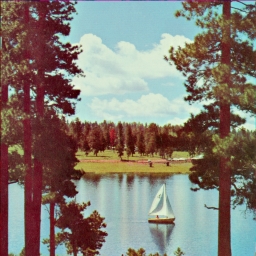}%
	\\
	\text{{\footnotesize (h)}}
\end{minipage}%
	\caption{The original color images with size of $256\times 256\times 3$. (a) {Airplane}, (b) {Baboon}, (c) {Barbara}, (d) {Facade}, (e) {House}, (f) {Lena}, (g) Peppers, (h) Sailboat.}
	\label{image_ori}
\end{figure}

We employ benchmark color images with size of $256\times 256 \times 3$ to test the recovering performance of algorithms under different SRs, and the original color images are demonstrated in Fig. \ref{image_ori}.

\begin{figure*}[htbp]
%%%%%%%%%%%%%%%%%%%%%%%%%%%%%%%%%%%%%%%%%%%%%%%%%%%%%%%%%%%%%%%%%%%%%%%%%%%%%%%%%%%%PSNR
	\centering	
	\begin{minipage}[f]{0.13\linewidth}
	\centering
	\includegraphics[width=1.5cm]{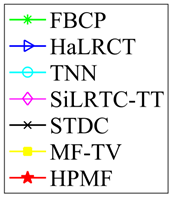}%
	\end{minipage}%
	\centering	
	\begin{minipage}[f]{0.21\linewidth}
	\centering
	\includegraphics[width=3.9cm]{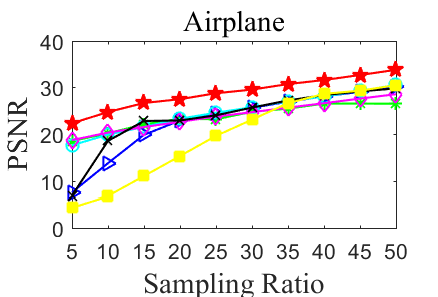}%
	\end{minipage}%
	\centering
	\begin{minipage}[f]{0.21\linewidth}
	\centering
	\includegraphics[width=3.9cm]{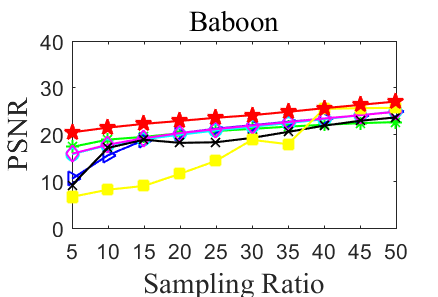}%
	\end{minipage}%
	\centering
	\begin{minipage}[f]{0.21\linewidth}
	\centering
	\includegraphics[width=3.9cm]{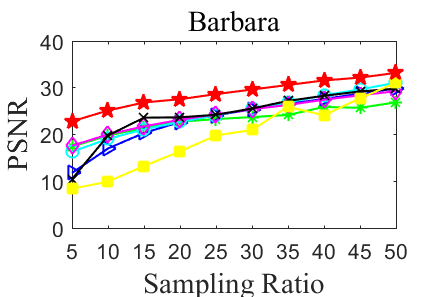}%	
	\end{minipage}%
	\centering
	\begin{minipage}[f]{0.21\linewidth}
	\centering
	\includegraphics[width=3.9cm]{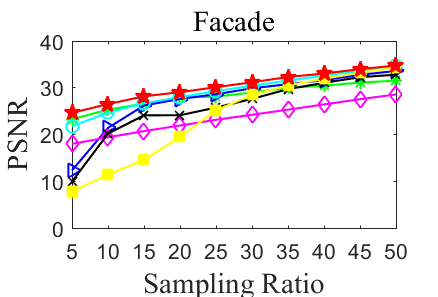}%
	\end{minipage}%
\\%%%%%%%%%%%%%%%%%%%%%%%%%%%%%%%%%%%%%%%%%%%%%%%%%%%%%%%%%%%%%
	\centering	
\begin{minipage}[f]{0.13\linewidth}
	\centering
	\includegraphics[width=1.5cm]{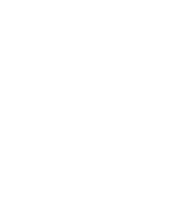}%
\end{minipage}%
	\centering	
\begin{minipage}[f]{0.21\linewidth}
	\centering
	\includegraphics[width=3.9cm]{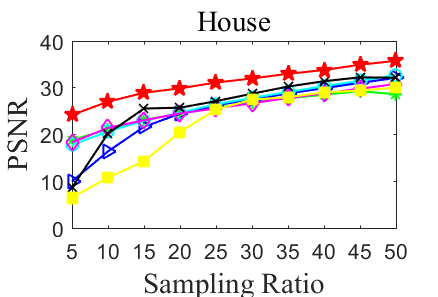}%
\end{minipage}%
\centering
\begin{minipage}[f]{0.21\linewidth}
	\centering
	\includegraphics[width=3.9cm]{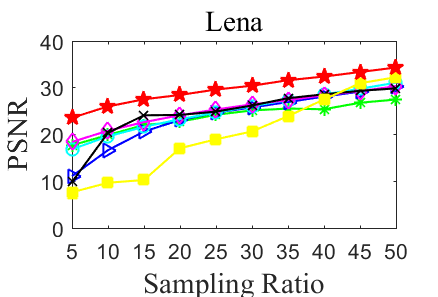}%
\end{minipage}%
\centering
\begin{minipage}[f]{0.21\linewidth}
	\centering
	\includegraphics[width=3.9cm]{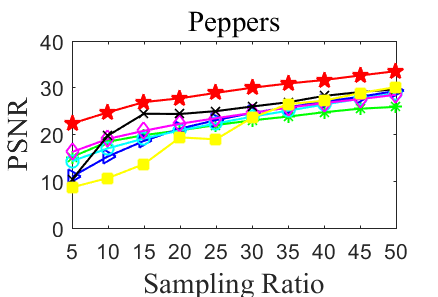}%	
\end{minipage}%
\centering
\begin{minipage}[f]{0.21\linewidth}
	\centering
	\includegraphics[width=3.9cm]{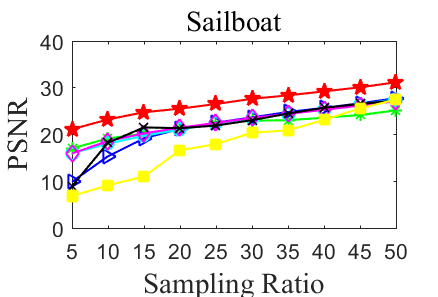}%
\end{minipage}%
\\
\vspace{0.2cm}
\text{{\footnotesize (a)}}
\vspace{0.4cm}
\\
%%%%%%%%%%%%%%%%%%%%%%%%%%%%%%%%%%%%%%%%%%%%%%%%%%%%%%%%%%%%%%%%%%%%%%%%%%%%%%%%%%%%RSE
	\centering	
\begin{minipage}[f]{0.13\linewidth}
	\centering
	\includegraphics[width=1.5cm]{legend.png}%
\end{minipage}%
\centering	
\begin{minipage}[f]{0.21\linewidth}
	\centering
	\includegraphics[width=3.9cm]{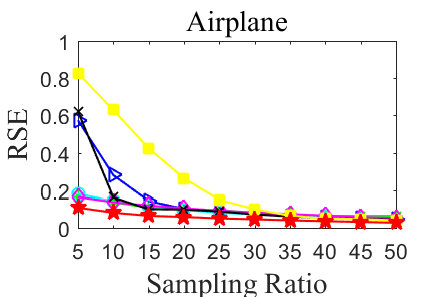}%
\end{minipage}%
\centering
\begin{minipage}[f]{0.21\linewidth}
	\centering
	\includegraphics[width=3.9cm]{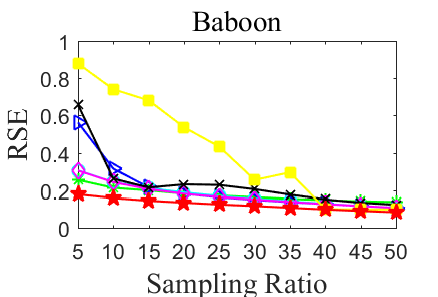}%
\end{minipage}%
\centering
\begin{minipage}[f]{0.21\linewidth}
	\centering
	\includegraphics[width=3.9cm]{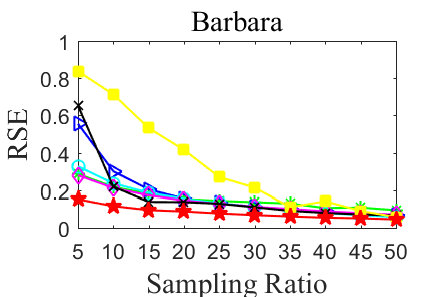}%	
\end{minipage}%
\centering
\begin{minipage}[f]{0.21\linewidth}
	\centering
	\includegraphics[width=3.9cm]{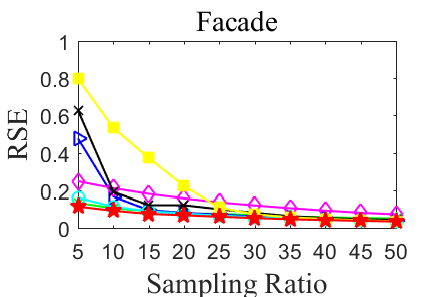}%
\end{minipage}%
\\%%%%%%%%%%%%%%%%%%%%%%%%%%%%%%%%%%%%%%%%%%%%%%%%%%%%%%%%%%%%%
\centering	
\begin{minipage}[f]{0.13\linewidth}
	\centering
	\includegraphics[width=1.9cm]{legend_.png}%
\end{minipage}%
\centering	
\begin{minipage}[f]{0.21\linewidth}
	\centering
	\includegraphics[width=3.9cm]{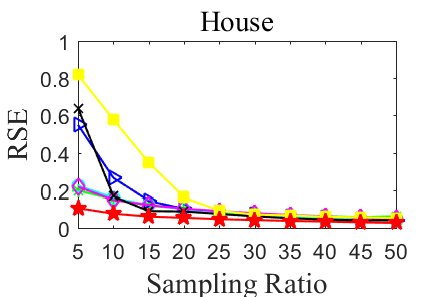}%
\end{minipage}%
\centering
\begin{minipage}[f]{0.21\linewidth}
	\centering
	\includegraphics[width=3.9cm]{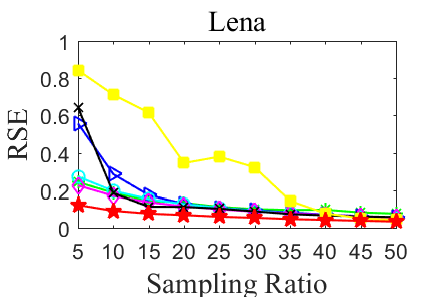}%
\end{minipage}%
\centering
\begin{minipage}[f]{0.21\linewidth}
	\centering
	\includegraphics[width=3.9cm]{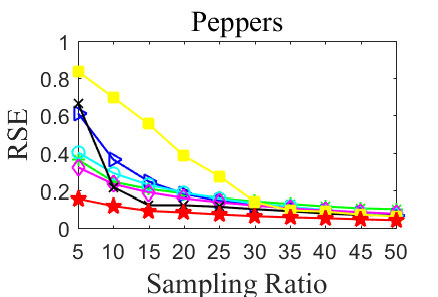}%	
\end{minipage}%
\centering
\begin{minipage}[f]{0.21\linewidth}
	\centering
	\includegraphics[width=3.9cm]{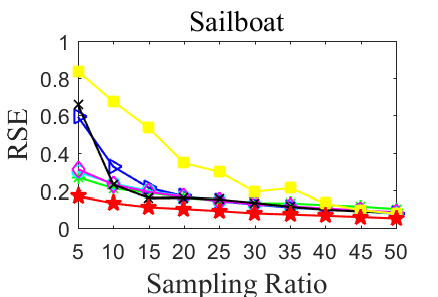}%
\end{minipage}%
\\
\vspace{0.2cm}
\text{{\footnotesize (b)}}
\vspace{0.4cm}
\\
%%%%%%%%%%%%%%%%%%%%%%%%%%%%%%%%%%%%%%%%%%%%%%%%%%%%%%%%%%%%%%%%%%%%%%%%%%%%%%%%%%%%SSIM
	\centering	
\begin{minipage}[f]{0.13\linewidth}
	\centering
	\includegraphics[width=1.9cm]{legend.png}%
\end{minipage}%
\centering	
\begin{minipage}[f]{0.21\linewidth}
	\centering
	\includegraphics[width=3.9cm]{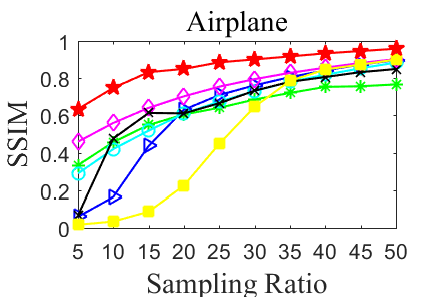}%
\end{minipage}%
\centering
\begin{minipage}[f]{0.21\linewidth}
	\centering
	\includegraphics[width=3.9cm]{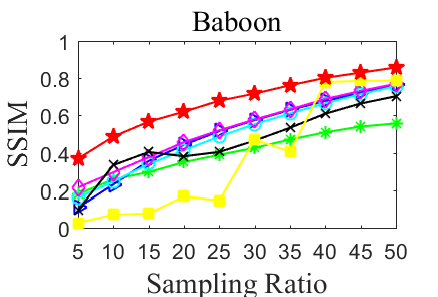}%
\end{minipage}%
\centering
\begin{minipage}[f]{0.21\linewidth}
	\centering
	\includegraphics[width=3.9cm]{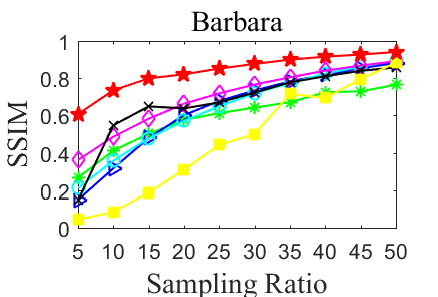}%	
\end{minipage}%
\centering
\begin{minipage}[f]{0.21\linewidth}
	\centering
	\includegraphics[width=3.9cm]{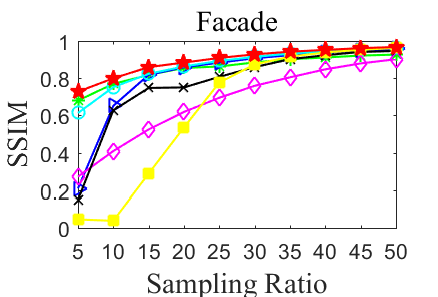}%
\end{minipage}%
\\%%%%%%%%%%%%%%%%%%%%%%%%%%%%%%%%%%%%%%%%%%%%%%%%%%%%%%%%%%%%%
\centering	
\begin{minipage}[f]{0.13\linewidth}
	\centering
	\includegraphics[width=1.9cm]{legend_.png}%
\end{minipage}%
\centering	
\begin{minipage}[f]{0.21\linewidth}
	\centering
	\includegraphics[width=3.9cm]{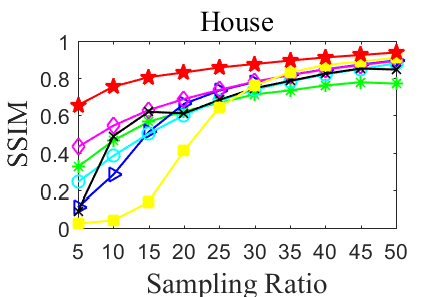}%
\end{minipage}%
\centering
\begin{minipage}[f]{0.21\linewidth}
	\centering
	\includegraphics[width=3.9cm]{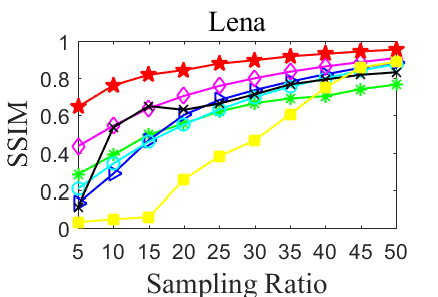}%
\end{minipage}%
\centering
\begin{minipage}[f]{0.21\linewidth}
	\centering
	\includegraphics[width=3.9cm]{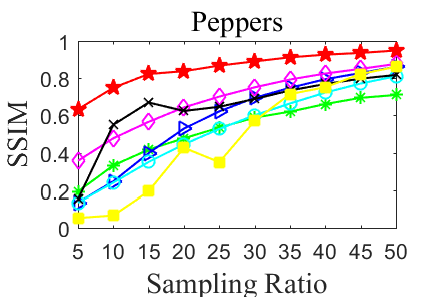}%	
\end{minipage}%
\centering
\begin{minipage}[f]{0.21\linewidth}
	\centering
	\includegraphics[width=3.9cm]{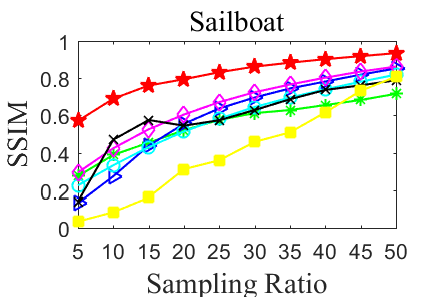}%
\end{minipage}%
\\
\vspace{0.2cm}
\text{{\footnotesize (c)}}
\vspace{0.2cm}
\\
\caption{Comparison of algorithms for color images with the SR from 5\% to 50\%. All compared algorithms are labeled in the legend. (a) Comparison for PSNR vs SR, (b) Comparison for RSE vs SR, (c) Comparison for SSIM vs SR.}
\label{CI_curves}
\end{figure*}

Fig. \ref{CI_curves} shows the results of tensor completion for benchmark color images when the SRs varies from 5\% to 50\%. In the figure, we compare the proposed HPMF algorithm and the state-of-the-art algorithms in terms of PSNR, RSE, and SSIM. We observe that, our HPMF method obtains the best values of PSNR, RSE, and SSIM, for all SRs, which demonstrates the HPMF effectively mines both the global low-rank structure and the local smoothness of benchmark images. Moreover, low-rank prior based methods (FBCP, TNN, SiLRTC-TT) exhibit good PSNR, RSE, and SSIM values, especially when the SR is below 20\%. HaLRTC only presents good performance with respect to PSNR, RSE, and SSIM, with the SR higher than 20\%. For smoothness regularized schemes, when the SR is above 30\%, STDC gains good PSNR and RSE values for most images, e.g., Barbara, House, and Peppers, which indicates the smoothness constraint may contribute more to the recovery quality given more observations. However, MF-TV depicts bad PSNR, RSE, and SSIM values to recover images, particularly when the SR is below 20\%. This phenomenon may explain that, the partial constraint smoothness for the mode-3 unfolding fails to sufficiently exploit the local patterns of images given limited observations.

\begin{figure}[H]
	\centering
	\begin{minipage}[f]{0.125\linewidth}
		\centering
		\includegraphics[width=1.05cm]{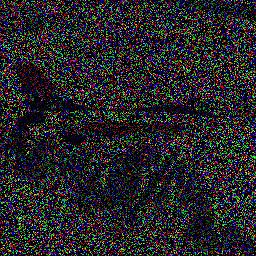}%
		\\
		%\text{{\footnotesize (a)}}
	\end{minipage}%
	\centering
\begin{minipage}[f]{0.125\linewidth}
	\centering
	\includegraphics[width=1.05cm]{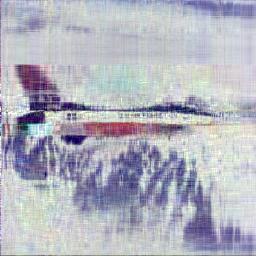}%
	\\
	%\text{{\footnotesize (b)}}
\end{minipage}%
	\centering
\begin{minipage}[f]{0.125\linewidth}
	\centering
	\includegraphics[width=1.05cm]{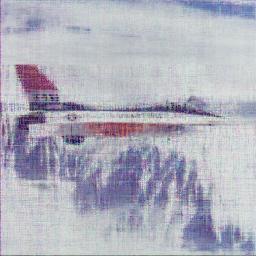}%
	\\
	%\text{{\footnotesize (c)}}
\end{minipage}%
	\centering
\begin{minipage}[f]{0.125\linewidth}
	\centering
	\includegraphics[width=1.05cm]{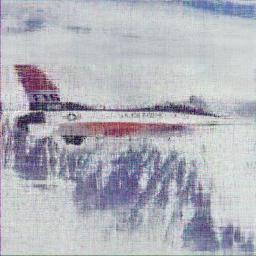}%
	\\
	%\text{{\footnotesize (d)}}
\end{minipage}%
	\centering
\begin{minipage}[f]{0.125\linewidth}
	\centering
	\includegraphics[width=1.05cm]{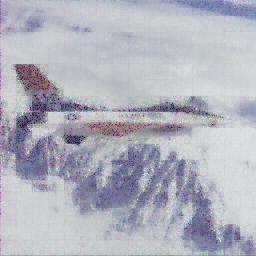}%
	\\
	%\text{{\footnotesize (e)}}
\end{minipage}%
	\centering
\begin{minipage}[f]{0.125\linewidth}
	\centering
	\includegraphics[width=1.05cm]{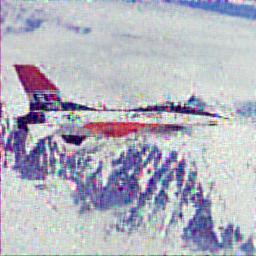}%
	\\
	%\text{{\footnotesize (f)}}
\end{minipage}%
	\centering
\begin{minipage}[f]{0.125\linewidth}
	\centering
	\includegraphics[width=1.05cm]{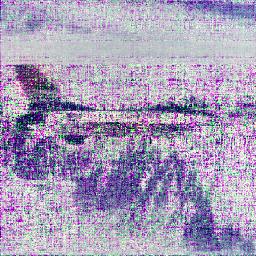}%
	\\
	%\text{{\footnotesize (g)}}
\end{minipage}%
	\centering
\begin{minipage}[f]{0.125\linewidth}
	\centering
	\includegraphics[width=1.05cm]{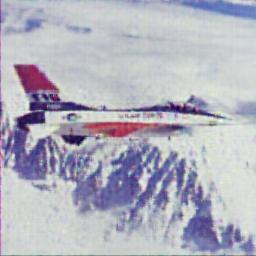}%
	\\
	%\text{{\footnotesize (h)}}
\end{minipage}%
\\
\vspace{1.5pt}%%%%%%%%%%%%%%%%%%%%%%%%%%%%%%%%%%%%%%%%%%%%%%%%%%%%%%%%%%%%%%%%%%%%%%%
	\centering
\begin{minipage}[f]{0.125\linewidth}
	\centering
	\includegraphics[width=1.05cm]{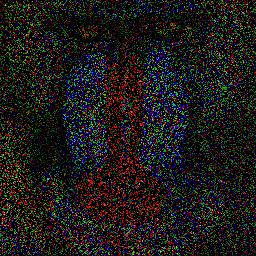}%
	\\
	%\text{{\footnotesize (a)}}
\end{minipage}%
\centering
\begin{minipage}[f]{0.125\linewidth}
	\centering
	\includegraphics[width=1.05cm]{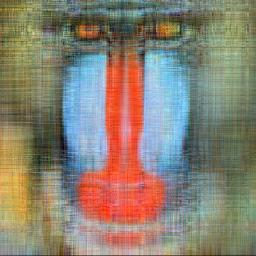}%
	\\
	%\text{{\footnotesize (b)}}
\end{minipage}%
\centering
\begin{minipage}[f]{0.125\linewidth}
	\centering
	\includegraphics[width=1.05cm]{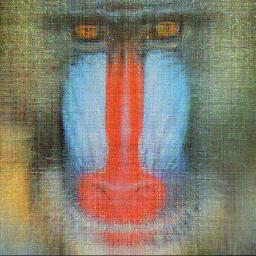}%
	\\
	%\text{{\footnotesize (c)}}
\end{minipage}%
\centering
\begin{minipage}[f]{0.125\linewidth}
	\centering
	\includegraphics[width=1.05cm]{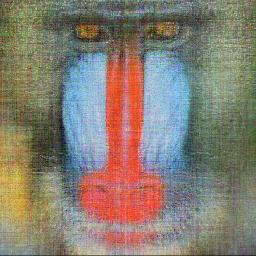}%
	\\
	%\text{{\footnotesize (d)}}
\end{minipage}%
\centering
\begin{minipage}[f]{0.125\linewidth}
	\centering
	\includegraphics[width=1.05cm]{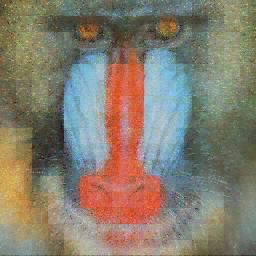}%
	\\
	%\text{{\footnotesize (e)}}
\end{minipage}%
\centering
\begin{minipage}[f]{0.125\linewidth}
	\centering
	\includegraphics[width=1.05cm]{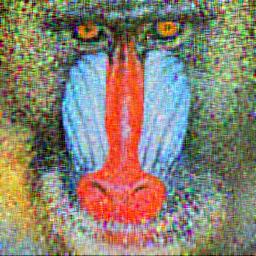}%
	\\
	%\text{{\footnotesize (f)}}
\end{minipage}%
\centering
\begin{minipage}[f]{0.125\linewidth}
	\centering
	\includegraphics[width=1.05cm]{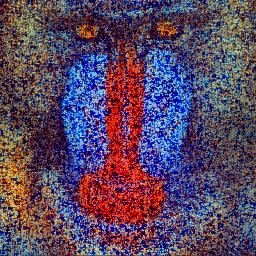}%
	\\
	%\text{{\footnotesize (g)}}
\end{minipage}%
\centering
\begin{minipage}[f]{0.125\linewidth}
	\centering
	\includegraphics[width=1.05cm]{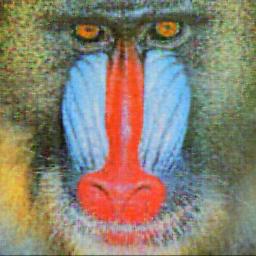}%
	\\
	%\text{{\footnotesize (h)}}
\end{minipage}%
\\
\vspace{1.5pt}%%%%%%%%%%%%%%%%%%%%%%%%%%%%%%%%%%%%%%%%%%%%%%%%%%%%
	\centering
\begin{minipage}[f]{0.125\linewidth}
	\centering
	\includegraphics[width=1.05cm]{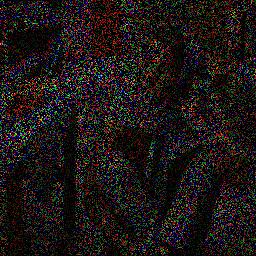}%
	\\
	%\text{{\footnotesize (a)}}
\end{minipage}%
\centering
\begin{minipage}[f]{0.125\linewidth}
	\centering
	\includegraphics[width=1.05cm]{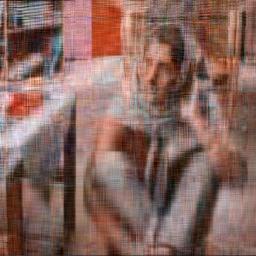}%
	\\
	%\text{{\footnotesize (b)}}
\end{minipage}%
\centering
\begin{minipage}[f]{0.125\linewidth}
	\centering
	\includegraphics[width=1.05cm]{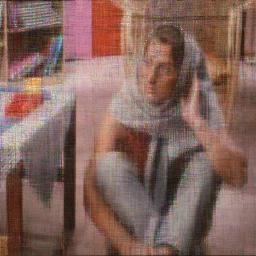}%
	\\
	%\text{{\footnotesize (c)}}
\end{minipage}%
\centering
\begin{minipage}[f]{0.125\linewidth}
	\centering
	\includegraphics[width=1.05cm]{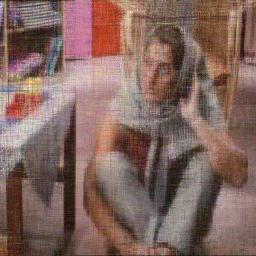}%
	\\
	%\text{{\footnotesize (d)}}
\end{minipage}%
\centering
\begin{minipage}[f]{0.125\linewidth}
	\centering
	\includegraphics[width=1.05cm]{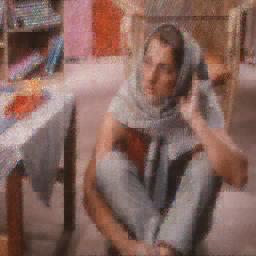}%
	\\
	%\text{{\footnotesize (e)}}
\end{minipage}%
\centering
\begin{minipage}[f]{0.125\linewidth}
	\centering
	\includegraphics[width=1.05cm]{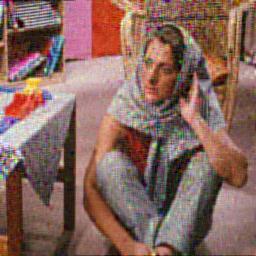}%
	\\
	%\text{{\footnotesize (f)}}
\end{minipage}%
\centering
\begin{minipage}[f]{0.125\linewidth}
	\centering
	\includegraphics[width=1.05cm]{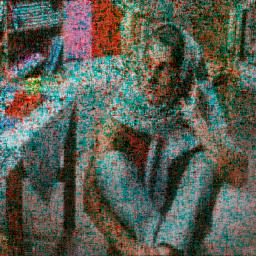}%
	\\
	%\text{{\footnotesize (g)}}
\end{minipage}%
\centering
\begin{minipage}[f]{0.125\linewidth}
	\centering
	\includegraphics[width=1.05cm]{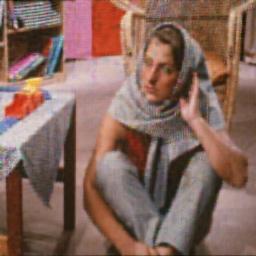}%
	\\
	%\text{{\footnotesize (h)}}
\end{minipage}%
\\
\vspace{1.5pt}%%%%%%%%%%%%%%%%%%%%%%%%%%%%%%%%%%%%%%%%%%%%%%%%%%%%
	\centering
\begin{minipage}[f]{0.125\linewidth}
	\centering
	\includegraphics[width=1.05cm]{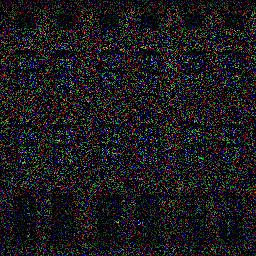}%
	\\
	%\text{{\footnotesize (a)}}
\end{minipage}%
\centering
\begin{minipage}[f]{0.125\linewidth}
	\centering
	\includegraphics[width=1.05cm]{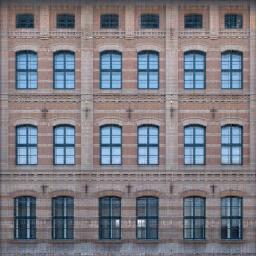}%
	\\
	%\text{{\footnotesize (b)}}
\end{minipage}%
\centering
\begin{minipage}[f]{0.125\linewidth}
	\centering
	\includegraphics[width=1.05cm]{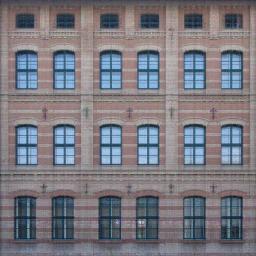}%
	\\
	%\text{{\footnotesize (c)}}
\end{minipage}%
\centering
\begin{minipage}[f]{0.125\linewidth}
	\centering
	\includegraphics[width=1.05cm]{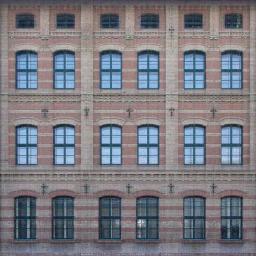}%
	\\
	%\text{{\footnotesize (d)}}
\end{minipage}%
\centering
\begin{minipage}[f]{0.125\linewidth}
	\centering
	\includegraphics[width=1.05cm]{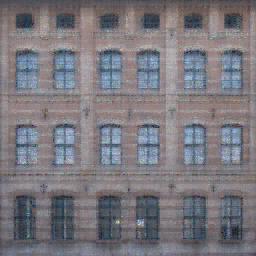}%
	\\
	%\text{{\footnotesize (e)}}
\end{minipage}%
\centering
\begin{minipage}[f]{0.125\linewidth}
	\centering
	\includegraphics[width=1.05cm]{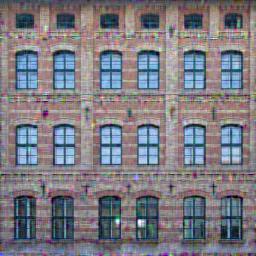}%
	\\
	%\text{{\footnotesize (f)}}
\end{minipage}%
\centering
\begin{minipage}[f]{0.125\linewidth}
	\centering
	\includegraphics[width=1.05cm]{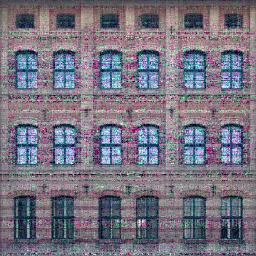}%
	\\
	%\text{{\footnotesize (g)}}
\end{minipage}%
\centering
\begin{minipage}[f]{0.125\linewidth}
	\centering
	\includegraphics[width=1.05cm]{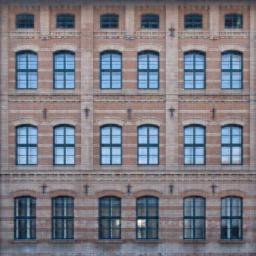}%
	\\
	%\text{{\footnotesize (h)}}
\end{minipage}%
\\
\vspace{1.5pt}%%%%%%%%%%%%%%%%%%%%%%%%%%%%%%%%%%%%%%%%%%%%%%%%%%%%
	\centering
\begin{minipage}[f]{0.125\linewidth}
	\centering
	\includegraphics[width=1.05cm]{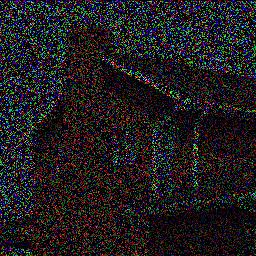}%
	\\
	%\text{{\footnotesize (a)}}
\end{minipage}%
\centering
\begin{minipage}[f]{0.125\linewidth}
	\centering
	\includegraphics[width=1.05cm]{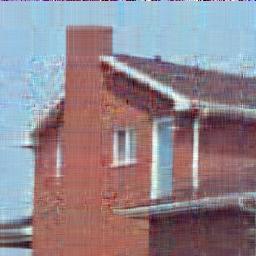}%
	\\
	%\text{{\footnotesize (b)}}
\end{minipage}%
\centering
\begin{minipage}[f]{0.125\linewidth}
	\centering
	\includegraphics[width=1.05cm]{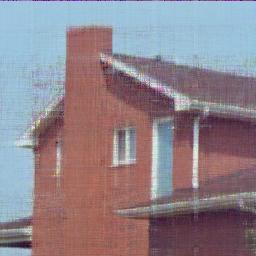}%
	\\
	%\text{{\footnotesize (c)}}
\end{minipage}%
\centering
\begin{minipage}[f]{0.125\linewidth}
	\centering
	\includegraphics[width=1.05cm]{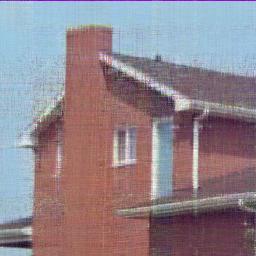}%
	\\
	%\text{{\footnotesize (d)}}
\end{minipage}%
\centering
\begin{minipage}[f]{0.125\linewidth}
	\centering
	\includegraphics[width=1.05cm]{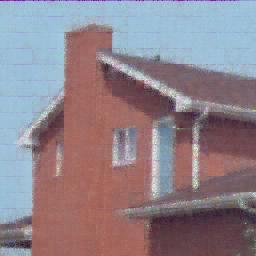}%
	\\
	%\text{{\footnotesize (e)}}
\end{minipage}%
\centering
\begin{minipage}[f]{0.125\linewidth}
	\centering
	\includegraphics[width=1.05cm]{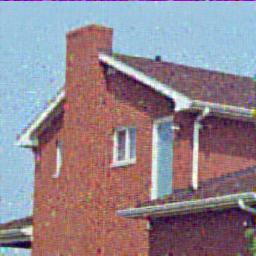}%
	\\
	%\text{{\footnotesize (f)}}
\end{minipage}%
\centering
\begin{minipage}[f]{0.125\linewidth}
	\centering
	\includegraphics[width=1.05cm]{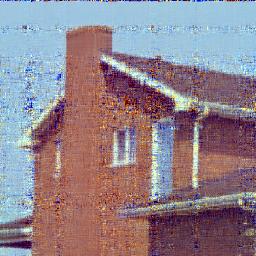}%
	\\
	%\text{{\footnotesize (g)}}
\end{minipage}%
\centering
\begin{minipage}[f]{0.125\linewidth}
	\centering
	\includegraphics[width=1.05cm]{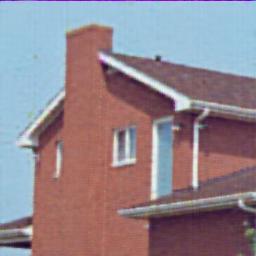}%
	\\
	%\text{{\footnotesize (h)}}
\end{minipage}%
\\
\vspace{1.5pt}%%%%%%%%%%%%%%%%%%%%%%%%%%%%%%%%%%%%%%%%%%%%%%%%%%%%
	\centering
\begin{minipage}[f]{0.125\linewidth}
	\centering
	\includegraphics[width=1.05cm]{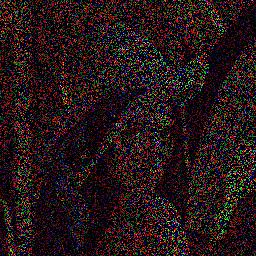}%
	\\
	%\text{{\footnotesize (a)}}
\end{minipage}%
\centering
\begin{minipage}[f]{0.125\linewidth}
	\centering
	\includegraphics[width=1.05cm]{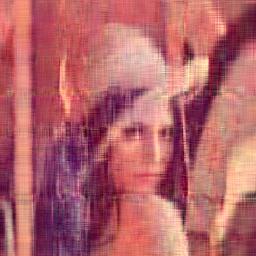}%
	\\
	%\text{{\footnotesize (b)}}
\end{minipage}%
\centering
\begin{minipage}[f]{0.125\linewidth}
	\centering
	\includegraphics[width=1.05cm]{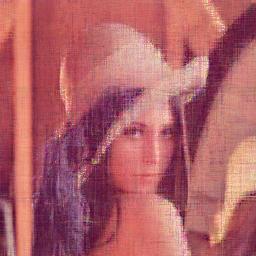}%
	\\
	%\text{{\footnotesize (c)}}
\end{minipage}%
\centering
\begin{minipage}[f]{0.125\linewidth}
	\centering
	\includegraphics[width=1.05cm]{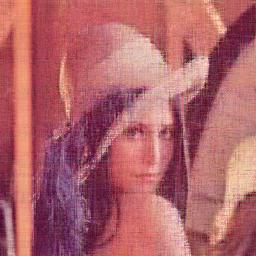}%
	\\
	%\text{{\footnotesize (d)}}
\end{minipage}%
\centering
\begin{minipage}[f]{0.125\linewidth}
	\centering
	\includegraphics[width=1.05cm]{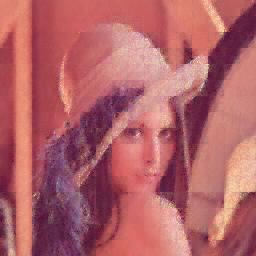}%
	\\
	%\text{{\footnotesize (e)}}
\end{minipage}%
\centering
\begin{minipage}[f]{0.125\linewidth}
	\centering
	\includegraphics[width=1.05cm]{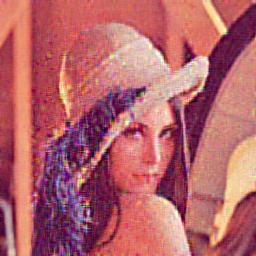}%
	\\
	%\text{{\footnotesize (f)}}
\end{minipage}%
\centering
\begin{minipage}[f]{0.125\linewidth}
	\centering
	\includegraphics[width=1.05cm]{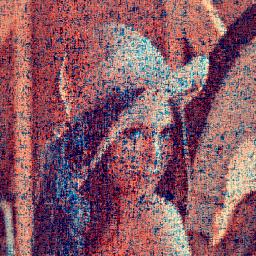}%
	\\
	%\text{{\footnotesize (g)}}
\end{minipage}%
\centering
\begin{minipage}[f]{0.125\linewidth}
	\centering
	\includegraphics[width=1.05cm]{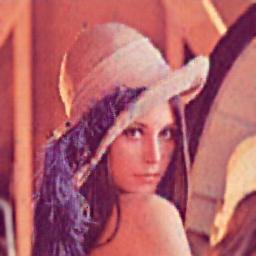}%
	\\
	%\text{{\footnotesize (h)}}
\end{minipage}%
\\
\vspace{1.5pt}%%%%%%%%%%%%%%%%%%%%%%%%%%%%%%%%%%%%%%%%%%%%%%%%%%%%
	\centering
\begin{minipage}[f]{0.125\linewidth}
	\centering
	\includegraphics[width=1.05cm]{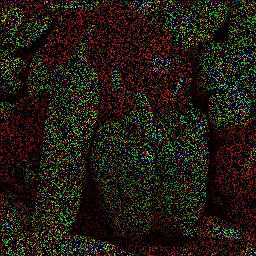}%
	\\
	%\text{{\footnotesize (a)}}
\end{minipage}%
\centering
\begin{minipage}[f]{0.125\linewidth}
	\centering
	\includegraphics[width=1.05cm]{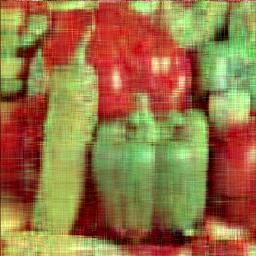}%
	\\
	%\text{{\footnotesize (b)}}
\end{minipage}%
\centering
\begin{minipage}[f]{0.125\linewidth}
	\centering
	\includegraphics[width=1.05cm]{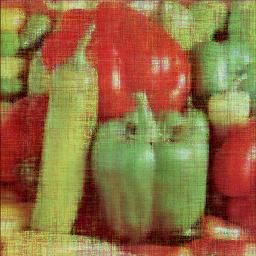}%
	\\
	%\text{{\footnotesize (c)}}
\end{minipage}%
\centering
\begin{minipage}[f]{0.125\linewidth}
	\centering
	\includegraphics[width=1.05cm]{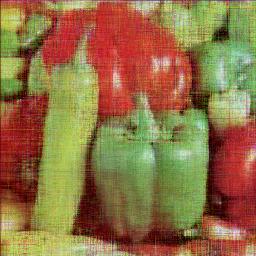}%
	\\
	%\text{{\footnotesize (d)}}
\end{minipage}%
\centering
\begin{minipage}[f]{0.125\linewidth}
	\centering
	\includegraphics[width=1.05cm]{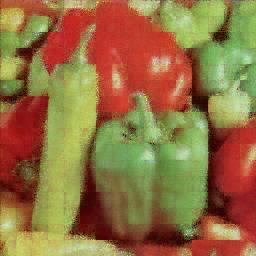}%
	\\
	%\text{{\footnotesize (e)}}
\end{minipage}%
\centering
\begin{minipage}[f]{0.125\linewidth}
	\centering
	\includegraphics[width=1.05cm]{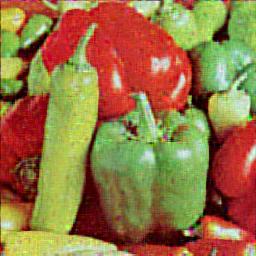}%
	\\
	%\text{{\footnotesize (f)}}
\end{minipage}%
\centering
\begin{minipage}[f]{0.125\linewidth}
	\centering
	\includegraphics[width=1.05cm]{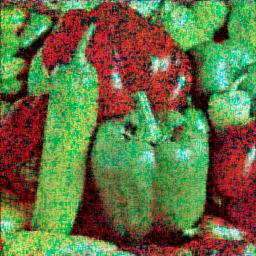}%
	\\
	%\text{{\footnotesize (g)}}
\end{minipage}%
\centering
\begin{minipage}[f]{0.125\linewidth}
	\centering
	\includegraphics[width=1.05cm]{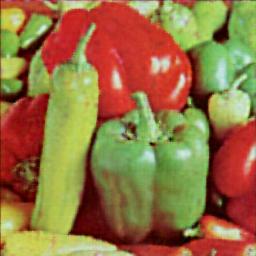}%
	\\
	%\text{{\footnotesize (h)}}
\end{minipage}%
\\
\vspace{1.5pt}%%%%%%%%%%%%%%%%%%%%%%%%%%%%%%%%%%%%%%%%%%%%%%%%%%%%
	\centering
\begin{minipage}[f]{0.125\linewidth}
	\centering
	\includegraphics[width=1.05cm]{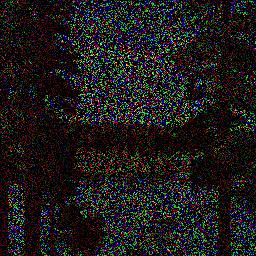}%
	\\
	\text{{\footnotesize (a)}}
\end{minipage}%
\centering
\begin{minipage}[f]{0.125\linewidth}
	\centering
	\includegraphics[width=1.05cm]{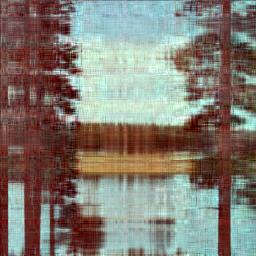}%
	\\
	\text{{\footnotesize (b)}}
\end{minipage}%
\centering
\begin{minipage}[f]{0.125\linewidth}
	\centering
	\includegraphics[width=1.05cm]{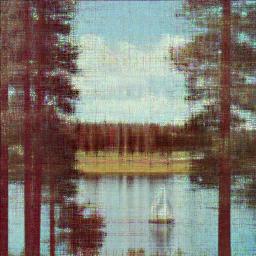}%
	\\
	\text{{\footnotesize (c)}}
\end{minipage}%
\centering
\begin{minipage}[f]{0.125\linewidth}
	\centering
	\includegraphics[width=1.05cm]{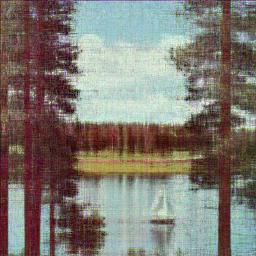}%
	\\
	\text{{\footnotesize (d)}}
\end{minipage}%
\centering
\begin{minipage}[f]{0.125\linewidth}
	\centering
	\includegraphics[width=1.05cm]{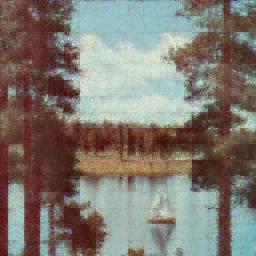}%
	\\
	\text{{\footnotesize (e)}}
\end{minipage}%
\centering
\begin{minipage}[f]{0.125\linewidth}
	\centering
	\includegraphics[width=1.05cm]{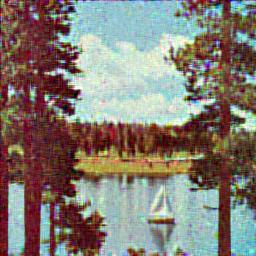}%
	\\
	\text{{\footnotesize (f)}}
\end{minipage}%
\centering
\begin{minipage}[f]{0.125\linewidth}
	\centering
	\includegraphics[width=1.05cm]{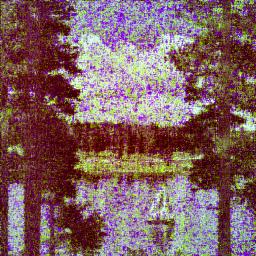}%
	\\
	\text{{\footnotesize (g)}}
\end{minipage}%
\centering
\begin{minipage}[f]{0.125\linewidth}
	\centering
	\includegraphics[width=1.05cm]{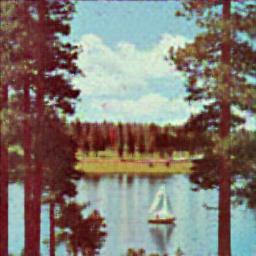}%
	\\
	\text{{\footnotesize (h)}}
\end{minipage}%
\\
\vspace{1.5pt}%%%%%%%%%%%%%%%%%%%%%%%%%%%%%%%%%%%%%%%%%%%%%%%%%%%%
	\caption{The observations and the recovered images via algorithms with the 20\% SR. From left to right: (a) {Observations}, (b) {FBCP}, (c) {HaLRTC}, (d) {TNN}, (e) {SiLRTC-TT}, (f) {STDC}, (g) MF-TV, (h) HPMF.}
	\label{CI_VC}
\end{figure}

In Fig. \ref{CI_VC}, we present the visual illustrations of recovered benchmark color images via compared algorithms when the SR is 20\%, to give the intuitive comparison. From this figure, we can see that, the proposed HPMF produces the best visual quality among all recovery approaches. Low-rank prior based methods (FBCP, HaLRTC, TNN, SiLRTC-TT) manifest good recovery of image outline, but they provide blurry image details. For the reconstructed images by STDC, they fail to show fine image textures. Besides, STDC tends to produce over-smoothed results, subsequently eliminating detail recovery of images. The recovered images via MF-TV is difficult to identify the image features, which further validates that, exclusive smooth constraint on mode-3 unfolding is not enough to gain necessary image details. Interestingly, for the images with obvious global low-rank structure, e.g., the Facade image, all algorithms show good performance of visual quality. In this scenario, the smooth constraint may be redundant, since it takes longer processing time to recover images. In summary, our HPMF method obtains clearer image details, sharper features, and finer textures, than compared schemes.

\begin{figure}[htbp]
	\centering
	\includegraphics[width=9.6cm]{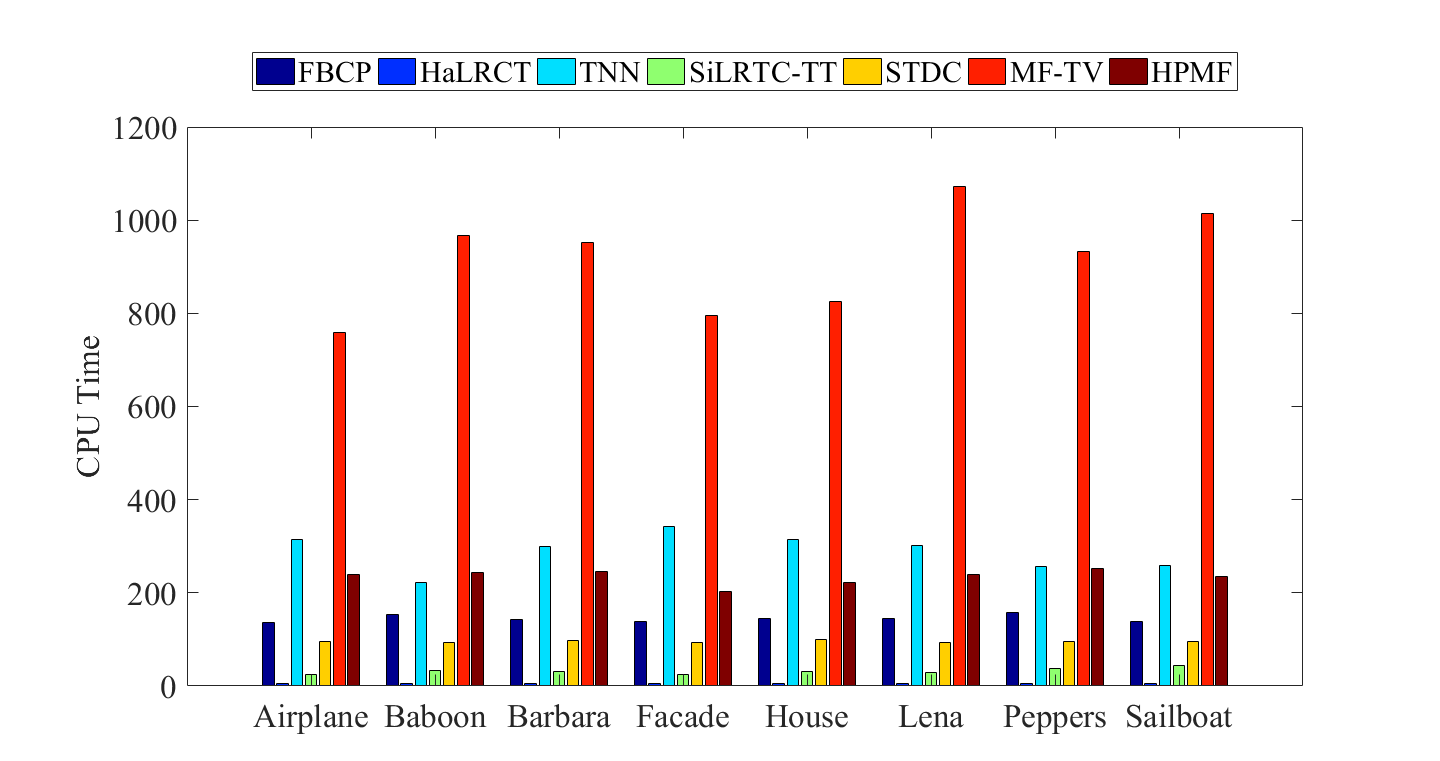}%
	\caption{The average CPU time (in seconds) for recovering color images via algorithms. All compared algorithms are labeled in the legend.}
	\label{CI_time}
\end{figure}

We provide the average CPU time (in seconds) to process benchmark color images with SRs varying from 5\% to 50\% in Fig. \ref{CI_time}. As illuminated in the figure, HaLRTC takes less CPU time than compared algorithms. Our HPMF costs moderate CPU time, which is less than TNN and MF-TV to recover most color images, except for the Baboon image. Due to the calculation of isotropic TV, MF-TV takes much longer processing time than other algorithms.

\subsection{Color Image with Non-random Sampling}
\begin{figure}[htbp]
	\centering
	\begin{minipage}[f]{0.33\linewidth}
		\centering
	\includegraphics[width=2.85cm]{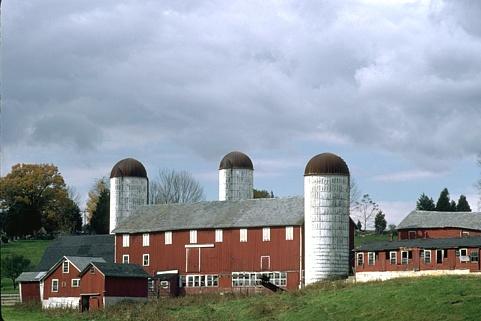}%
		\\
		\text{{\footnotesize (a)}}
	\end{minipage}%
	\begin{minipage}[f]{0.33\linewidth}
		\centering
		\includegraphics[width=2.85cm]{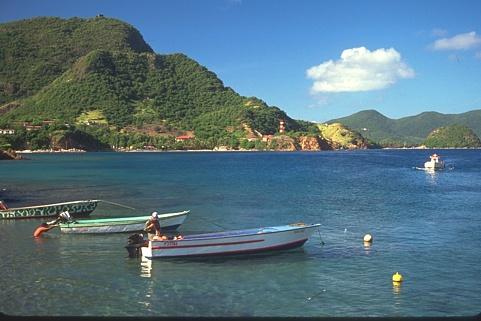}%
		\\
		\text{{\footnotesize (b)}}
	\end{minipage}%
	\begin{minipage}[f]{0.33\linewidth}
		\centering
		\includegraphics[width=2.85cm]{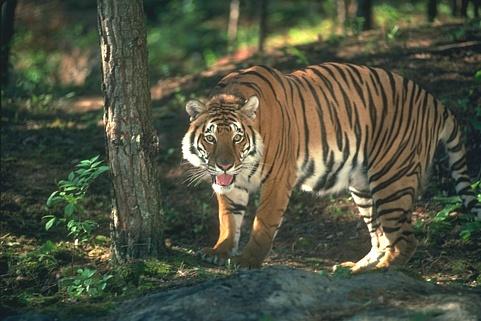}%
		\\
		\text{{\footnotesize (c)}}
	\end{minipage}%
\\
	\begin{minipage}[f]{0.33\linewidth}
		\centering
		\includegraphics[width=2.85cm]{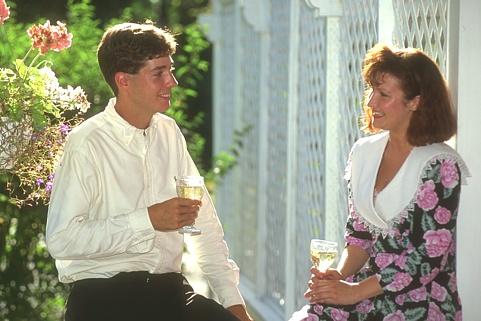}%
		\\
		\text{{\footnotesize (d)}}
	\end{minipage}%
	\begin{minipage}[f]{0.33\linewidth}
	\centering
	\includegraphics[width=2.85cm]{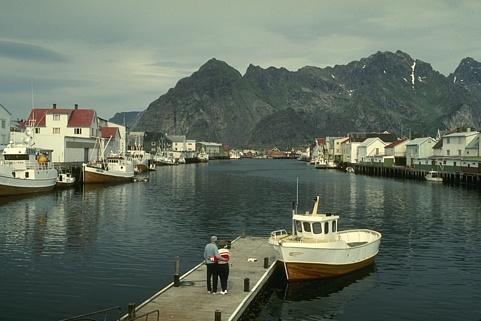}%
	\\
	\text{{\footnotesize (e)}}
\end{minipage}%
	\begin{minipage}[f]{0.33\linewidth}
	\centering
	\includegraphics[width=2.85cm]{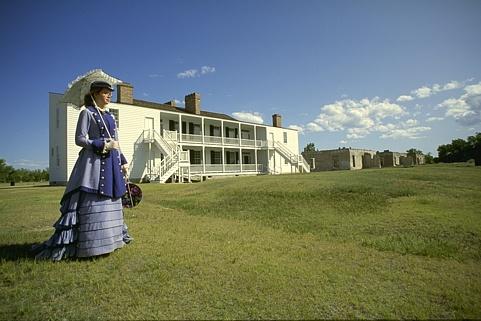}%
	\\
	\text{{\footnotesize (f)}}
\end{minipage}%
	\caption{The original images with size of $321\times 481\times 3$. (a) {Farm}, (b) {Coast}, (c) {Tiger}, (d) {People}, (e) {Port}, (f) {Lady}.}
	\label{image1_ori}
\end{figure}

This subsection conducts the experiments under various image masks. We utilize another group of color images$\footnote{https://www2.eecs.berkeley.edu/Research/Projects/CS/vision/bsds/.}$ with size of $321\times 481\times 3$ to exploit the recovering performance of algorithms, and Fig. \ref{image1_ori} illustrates the original color images. The settings are the same as in the color image with random sampling.

In Table \ref{maskTable}, we describe the recovery performance in terms of PSNR, RSE, SSIM, and average CPU time for color images with different image masks. The best PSNR, RSE, SSIM, and average CPU time are highlighted in boldface. It can be observed that, the proposed HPMF outperforms compared algorithms in terms of PSNR and RSE values for all masks. Furthermore, the SSIM values of HPMF are larger than compared state-of-the-art methods in most cases. All low-rank prior based schemes (FBCP, HaLRTC, TNN, SiLRTC-TT) remove masks efficiently. In detail, TNN and SiLRTC-TT produce good PSNR, SSIM, and RSE values for all masks. In all experiments, FBCP gives moderate performance in terms of PSNR, SSIM, and RSE. Although HaLRTC shows ordinary PSNR, RSE, and SSIM values, it has the second best performance for removing the Alphabet mask. STDC fails to remove masks in the experiments, and MF-TV only demonstrates good PSNR, RSE, and SSIM values to remove the Circle mask. On average, the PSNR improvement of our HPMF method over the second best algorithm (TNN) is 1.549 dB. The average RSE decline of the HPMF over the second best scheme (TNN) is 0.012, and the average SSIM improvement of the HPMF over the second best algorithm (SiLRTC-TT) is 0.017. Therefore, it can be concluded that, the proposed HPMF for image completion with non-random sampling yields superior performance. As for the average CPU time to process images with masks, SiLRTC-TT costs the least CPU time among compared algorithms. Our HPMF takes moderate time in this scenario.

\begin{table*}[htb]
	\centering
	\setlength{\tabcolsep}{0.3cm}
	\caption{Recovery Results via Compared Algorithms for Color Images with Non-random Sampling}
	\label{maskTable}
	\begin{tabular}{clllllllc}
		\toprule
		%Algorithms& &Grid-Coast & Scratch-Tiger & Text-Farm & Line-Port & Circle-People & Alphabet-Lady & CPU Time \\ 
		Algorithm&  & Grid & Scratch & Text & Line & Circle & Alphabet & CPU Time \\ 
		\midrule	
		\multirow{3}{*}{FBCP} & PSNR & 28.198 & 24.603 & 29.035 & 25.187 & 26.264 & 28.454 & \multirow{3}{*}{523.7} \\ 
		& RSE & 0.082 & 0.174 & 0.054 & 0.122 & 0.077 & 0.078 &  \\ 
		& SSIM & 0.806 & 0.758 & 0.873 & 0.771 & 0.801 & 0.828 &  \\ 
		\midrule
		\multirow{3}{*}{HaLRTC} & PSNR & 25.778 & 26.818 & 31.406 & 27.368 & 28.333 & 33.064 & \multirow{3}{*}{8.6} \\ 
		& RSE & 0.110 & 0.135 & 0.041 & 0.095 & 0.061 & 0.046 &  \\ 
		& SSIM & 0.830 & 0.878 & 0.948 & 0.840 & 0.910 & 0.961 &  \\ 
		\midrule
		\multirow{3}{*}{TNN} & PSNR & 33.456 & 26.766 & 31.404 & 30.494 & 30.476 & 32.989 & \multirow{3}{*}{702.9} \\ 
		& RSE & 0.045 & 0.136 & 0.041 & 0.066 & 0.048 & 0.046 &  \\ 
		& SSIM & 0.951 & 0.869 & 0.948 & 0.924 & 0.937 & 0.960 &  \\ 
		\midrule
		\multirow{3}{*}{SiLRTC-TT} & PSNR & 32.912 & 26.408 & 30.095 & 30.247 & 31.033 & 31.786 & \multirow{3}{*}{\textbf{3.8}} \\
		& RSE & 0.048 & 0.142 & 0.048 & 0.068 & 0.045 & 0.053 &  \\
		& SSIM & 0.953 & 0.869 & 0.935 & 0.921 & \textbf{0.956} & 0.957 &  \\ 
		\midrule
		\multirow{3}{*}{STDC} & PSNR & 12.927 & 9.895 & 6.340 & 11.257 & 9.633 & 10.373 & \multirow{3}{*}{203.2} \\
		& RSE & 0.487 & 0.958 & 0.745 & 0.608 & 0.526 & 0.625 &  \\ 
		& SSIM & 0.192 & 0.040 & 0.046 & 0.145 & 0.212 & 0.149 &  \\
		\midrule
		\multirow{3}{*}{MF-TV} & PSNR & 25.928 & 22.265 & 30.277 & 23.653 & 27.574 & 28.050 & \multirow{3}{*}{1284.6} \\ 
		& RSE & 0.107 & 0.228 & 0.047 & 0.146 & 0.067 & 0.082 &  \\
		& SSIM & 0.849 & 0.778 & 0.915 & 0.741 & 0.901 & 0.897 &  \\
		\midrule
		\multirow{3}{*}{HPMF} & PSNR & \textbf{34.191} & \textbf{29.458} & \textbf{33.774} & \textbf{31.914} & \textbf{31.609} & \textbf{33.936} & \multirow{3}{*}{738.9} \\
		& RSE & \textbf{0.041} & \textbf{0.099} & \textbf{0.032} & \textbf{0.056} & \textbf{0.042} & \textbf{0.041} &  \\ 
		& SSIM & \textbf{0.959} & \textbf{0.926} & \textbf{0.964} & \textbf{0.932} & 0.948 & \textbf{0.963} &  \\ 
		\bottomrule
	\end{tabular}
\end{table*}
 
\begin{figure*}[htbp]
	\centering
	\begin{minipage}[f]{0.123\linewidth}
		\centering
		\text{Grid Mask}
		\includegraphics[width=2.15cm]{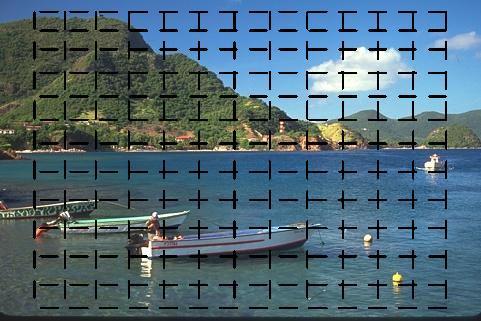}%
		\\
		%\text{{\footnotesize (a)}}
	\end{minipage}%
	\begin{minipage}[f]{0.123\linewidth}
		\centering
		\vspace{6.6pt}
		\includegraphics[width=2.15cm]{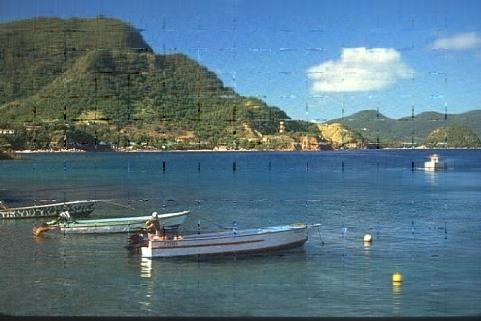}%
		\\
		%\text{{\footnotesize (b)}}
	\end{minipage}%
	\begin{minipage}[f]{0.123\linewidth}
		\centering
		\vspace{6.6pt}
		\includegraphics[width=2.15cm]{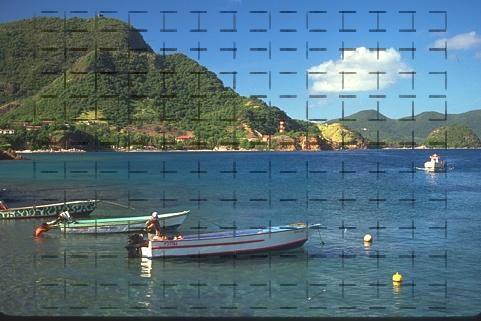}%
		\\
		%\text{{\footnotesize (c)}}
	\end{minipage}%
	\begin{minipage}[f]{0.123\linewidth}
		\centering
		\vspace{6.6pt}
		\includegraphics[width=2.15cm]{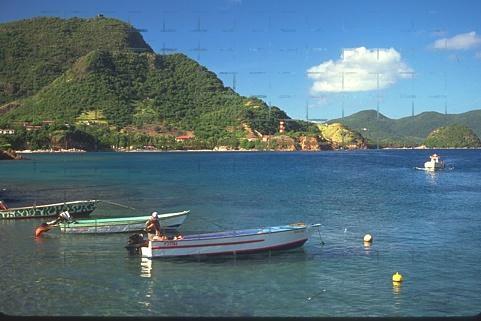}%
		\\
		%\text{{\footnotesize (d)}}
	\end{minipage}%
	\begin{minipage}[f]{0.123\linewidth}
		\centering
		\vspace{6.6pt}
		\includegraphics[width=2.15cm]{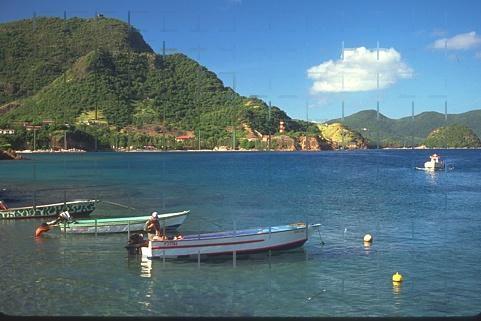}%
		\\
		%\text{{\footnotesize (e)}}
	\end{minipage}%
	\begin{minipage}[f]{0.123\linewidth}
		\centering
		\vspace{6.6pt}
		\includegraphics[width=2.15cm]{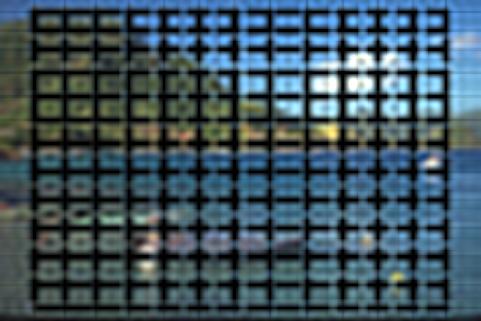}%
		\\
		%\text{{\footnotesize (f)}}
	\end{minipage}%
	\begin{minipage}[f]{0.123\linewidth}
	\centering
	\vspace{6.6pt}
	\includegraphics[width=2.15cm]{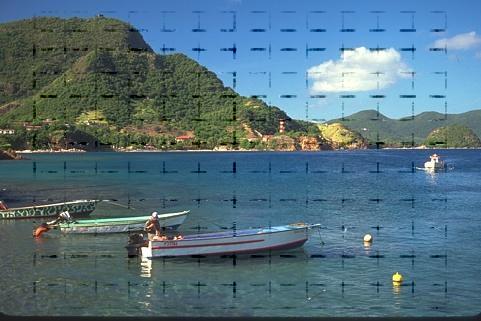}%
	\\
	%\text{{\footnotesize (g)}}
\end{minipage}%
	\begin{minipage}[f]{0.123\linewidth}
	\centering
	\vspace{6.6pt}
	\includegraphics[width=2.15cm]{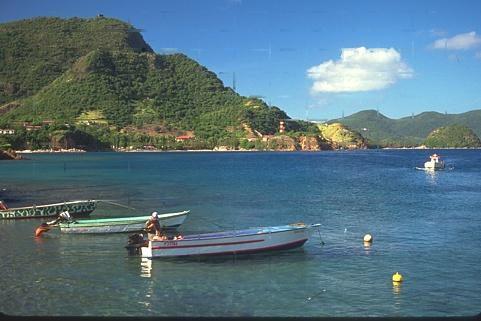}%
	\\
	%\text{{\footnotesize (h)}}
\end{minipage}%
\\
\vspace{1.5pt}%%%%%%%%%%%%%%%%%%%%%%%%%%%%%%%%%%%%%%%%%%%%%%%%%%%%
	\centering
\begin{minipage}[f]{0.123\linewidth}
	\centering
	\text{Scratch Mask}
	\includegraphics[width=2.15cm]{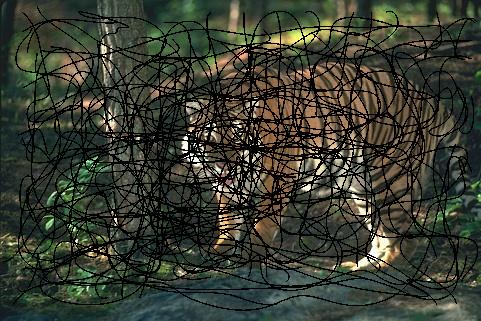}%
	\\
	%\text{{\footnotesize (a)}}
\end{minipage}%
\begin{minipage}[f]{0.123\linewidth}
	\centering
	\vspace{6.6pt}
	\includegraphics[width=2.15cm]{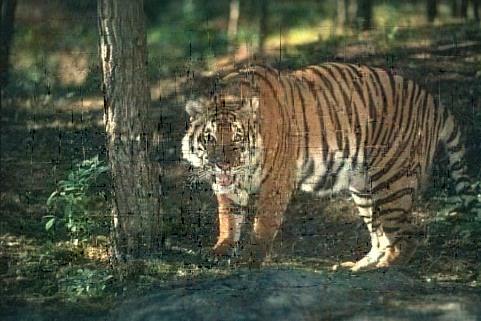}%
	\\
	%\text{{\footnotesize (b)}}
\end{minipage}%
\begin{minipage}[f]{0.123\linewidth}
	\centering
	\vspace{6.6pt}
	\includegraphics[width=2.15cm]{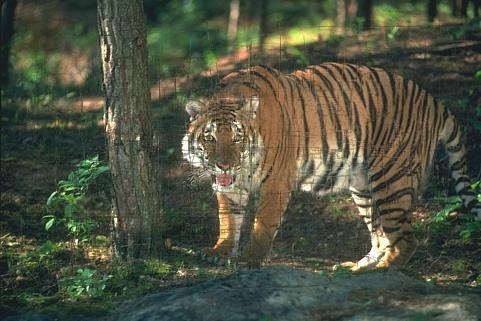}%
	\\
	%\text{{\footnotesize (c)}}
\end{minipage}%
\begin{minipage}[f]{0.123\linewidth}
	\centering
	\vspace{6.6pt}
	\includegraphics[width=2.15cm]{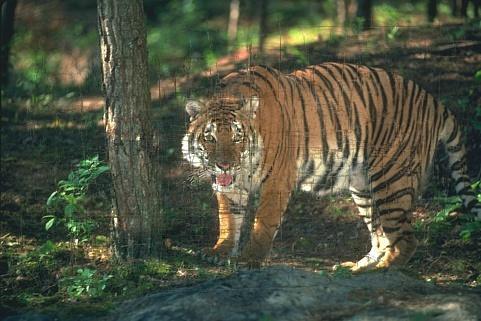}%
	\\
	%\text{{\footnotesize (d)}}
\end{minipage}%
\begin{minipage}[f]{0.123\linewidth}
	\centering
	\vspace{6.6pt}
	\includegraphics[width=2.15cm]{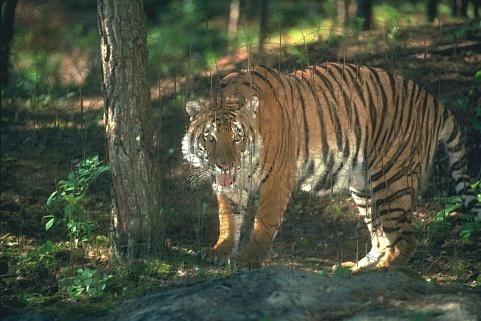}%
	\\
	%\text{{\footnotesize (e)}}
\end{minipage}%
\begin{minipage}[f]{0.123\linewidth}
	\centering
	\vspace{6.6pt}
	\includegraphics[width=2.15cm]{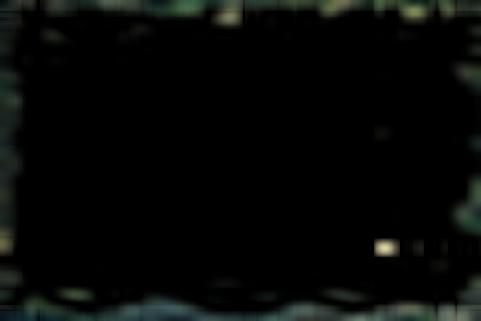}%
	\\
	%\text{{\footnotesize (f)}}
\end{minipage}%
\begin{minipage}[f]{0.123\linewidth}
	\centering
	\vspace{6.6pt}
	\includegraphics[width=2.15cm]{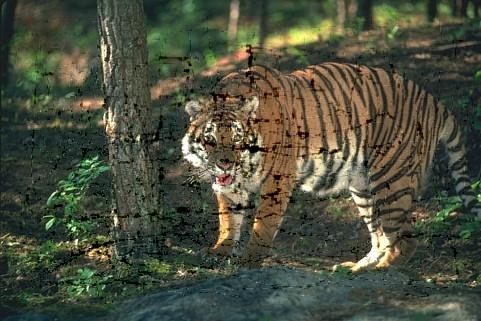}%
	\\
	%\text{{\footnotesize (g)}}
\end{minipage}%
\begin{minipage}[f]{0.123\linewidth}
	\centering
	\vspace{6.6pt}
	\includegraphics[width=2.15cm]{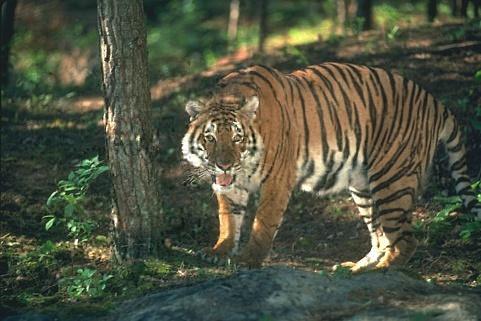}%
	\\
	%\text{{\footnotesize (h)}}
\end{minipage}%
\\
\vspace{1.5pt}%%%%%%%%%%%%%%%%%%%%%%%%%%%%%%%%%%%%%%%%%%%%%%%%%%%%
	\centering
\begin{minipage}[f]{0.123\linewidth}
	\centering
	\text{Text Mask}
	\includegraphics[width=2.15cm]{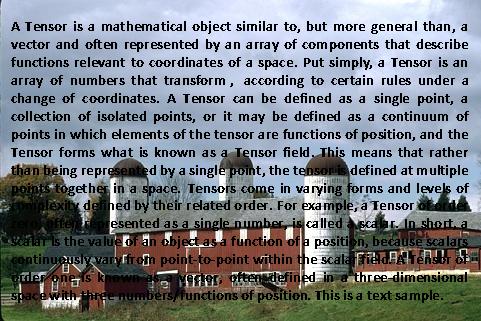}%
	\\
	%\text{{\footnotesize (a)}}
\end{minipage}%
\begin{minipage}[f]{0.123\linewidth}
	\centering
	\vspace{6.6pt}
	\includegraphics[width=2.15cm]{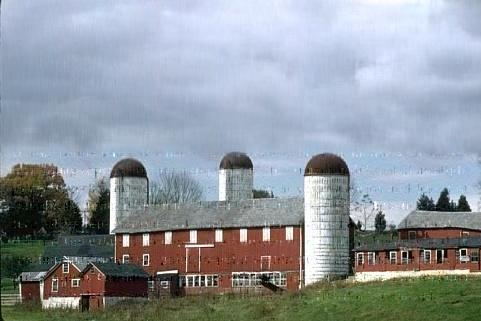}%
	\\
	%\text{{\footnotesize (b)}}
\end{minipage}%
\begin{minipage}[f]{0.123\linewidth}
	\centering
	\vspace{6.6pt}
	\includegraphics[width=2.15cm]{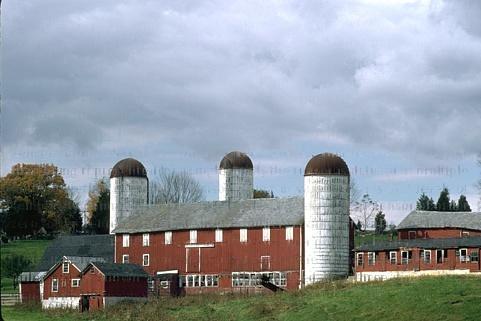}%
	\\
	%\text{{\footnotesize (c)}}
\end{minipage}%
\begin{minipage}[f]{0.123\linewidth}
	\centering
	\vspace{6.6pt}
	\includegraphics[width=2.15cm]{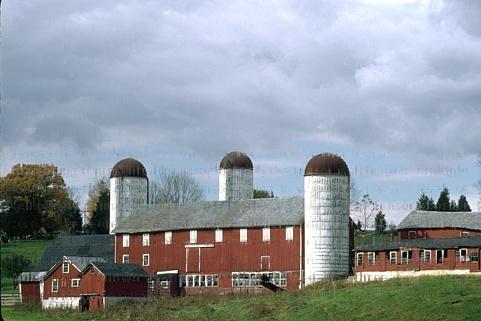}%
	\\
	%\text{{\footnotesize (d)}}
\end{minipage}%
\begin{minipage}[f]{0.123\linewidth}
	\centering
	\vspace{6.6pt}
	\includegraphics[width=2.15cm]{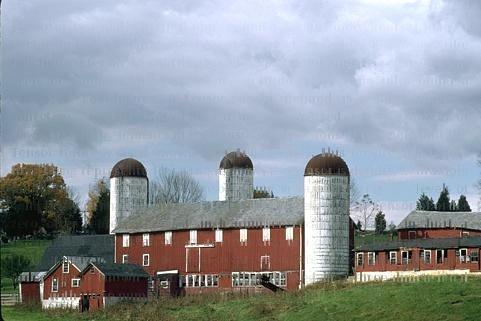}%
	\\
	%\text{{\footnotesize (e)}}
\end{minipage}%
\begin{minipage}[f]{0.123\linewidth}
	\centering
	\vspace{6.6pt}
	\includegraphics[width=2.15cm]{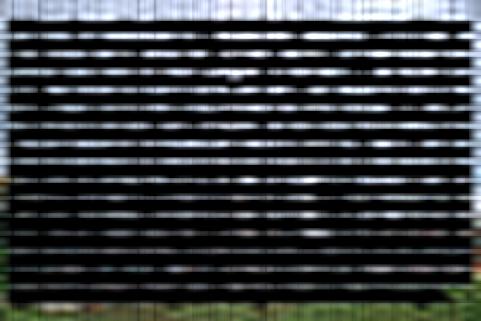}%
	\\
	%\text{{\footnotesize (f)}}
\end{minipage}%
\begin{minipage}[f]{0.123\linewidth}
	\centering
	\vspace{6.6pt}
	\includegraphics[width=2.15cm]{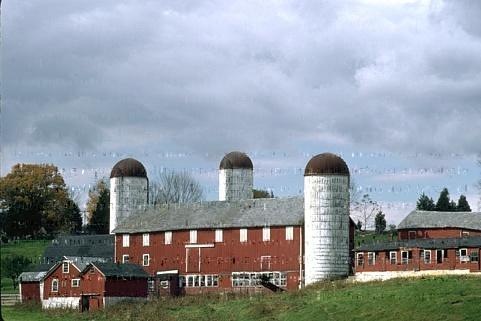}%
	\\
	%\text{{\footnotesize (g)}}
\end{minipage}%
\begin{minipage}[f]{0.123\linewidth}
	\centering
	\vspace{6.6pt}
	\includegraphics[width=2.15cm]{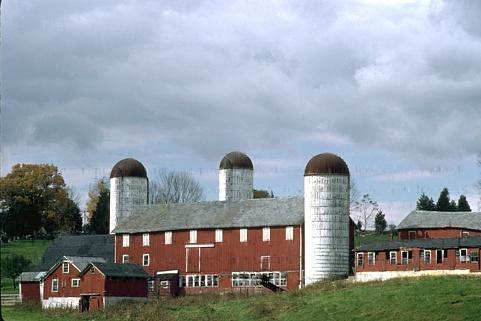}%
	\\
	%\text{{\footnotesize (h)}}
\end{minipage}%
\\
\vspace{1.5pt}%%%%%%%%%%%%%%%%%%%%%%%%%%%%%%%%%%%%%%%%%%%%%%%%%%%%
	\centering
\begin{minipage}[f]{0.123\linewidth}
	\centering
	\text{Line Mask}
	\includegraphics[width=2.15cm]{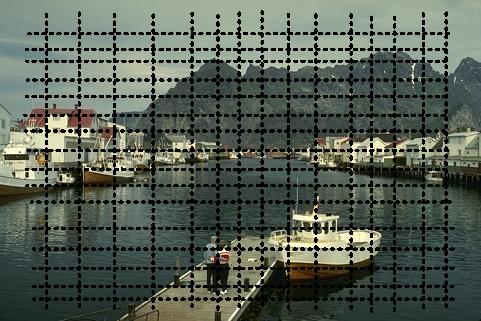}%
	\\
	%\text{{\footnotesize (a)}}
\end{minipage}%
\begin{minipage}[f]{0.123\linewidth}
	\centering
	\vspace{6.6pt}
	\includegraphics[width=2.15cm]{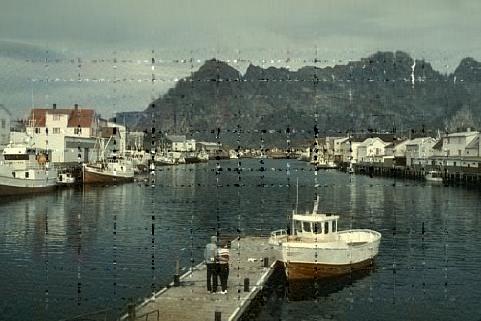}%
	\\
	%\text{{\footnotesize (b)}}
\end{minipage}%
\begin{minipage}[f]{0.123\linewidth}
	\centering
	\vspace{6.6pt}
	\includegraphics[width=2.15cm]{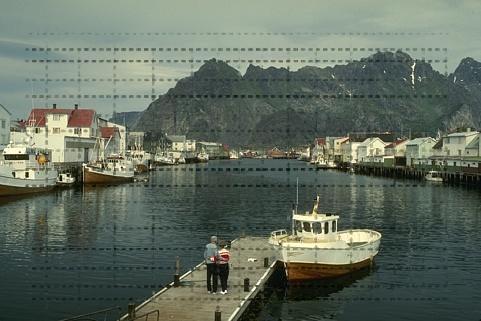}%
	\\
	%\text{{\footnotesize (c)}}
\end{minipage}%
\begin{minipage}[f]{0.123\linewidth}
	\centering
	\vspace{6.6pt}
	\includegraphics[width=2.15cm]{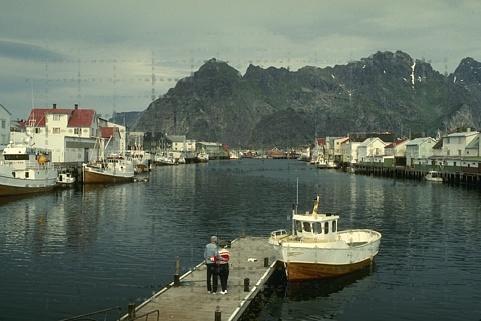}%
	\\
	%\text{{\footnotesize (d)}}
\end{minipage}%
\begin{minipage}[f]{0.123\linewidth}
	\centering
	\vspace{6.6pt}
	\includegraphics[width=2.15cm]{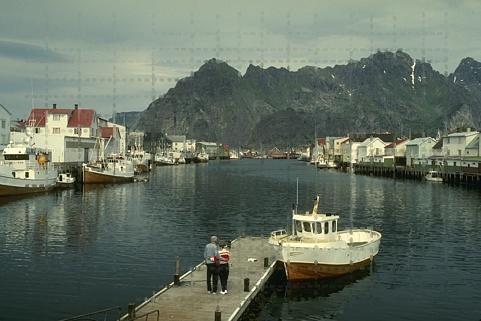}%
	\\
	%\text{{\footnotesize (e)}}
\end{minipage}%
\begin{minipage}[f]{0.123\linewidth}
	\centering
	\vspace{6.6pt}
	\includegraphics[width=2.15cm]{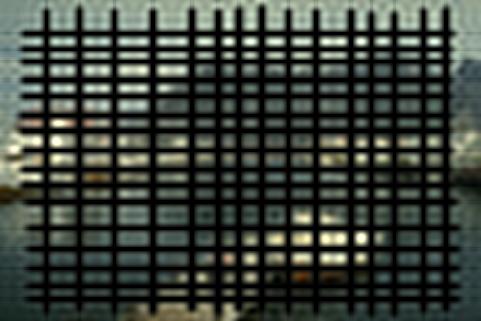}%
	\\
	%\text{{\footnotesize (f)}}
\end{minipage}%
\begin{minipage}[f]{0.123\linewidth}
	\centering
	\vspace{6.6pt}
	\includegraphics[width=2.15cm]{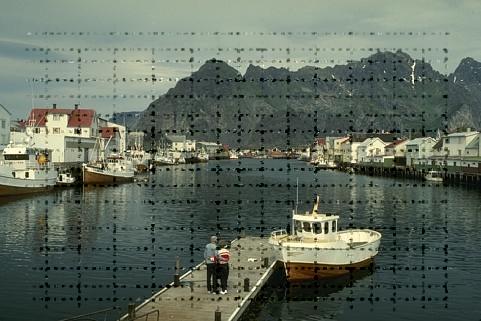}%
	\\
	%\text{{\footnotesize (g)}}
\end{minipage}%
\begin{minipage}[f]{0.123\linewidth}
	\centering
	\vspace{6.6pt}
	\includegraphics[width=2.15cm]{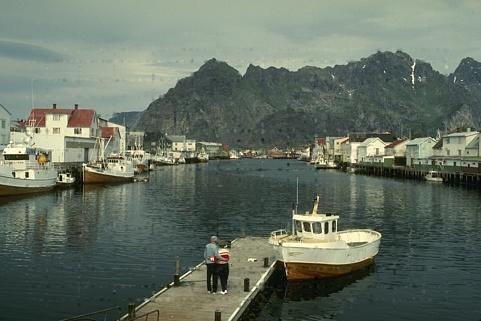}%
	\\
	%\text{{\footnotesize (h)}}
\end{minipage}%
\\
\vspace{1.5pt}%%%%%%%%%%%%%%%%%%%%%%%%%%%%%%%%%%%%%%%%%%%%%%%%%%%%
	\centering
\begin{minipage}[f]{0.123\linewidth}
	\centering
	\text{Circle Mask}
	\includegraphics[width=2.15cm]{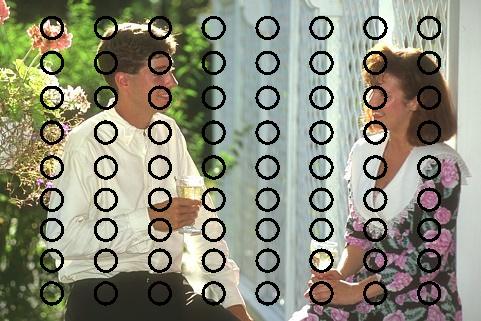}%
	\\
	%\text{{\footnotesize (a)}}
\end{minipage}%
\begin{minipage}[f]{0.123\linewidth}
	\centering
	\vspace{6.6pt}
	\includegraphics[width=2.15cm]{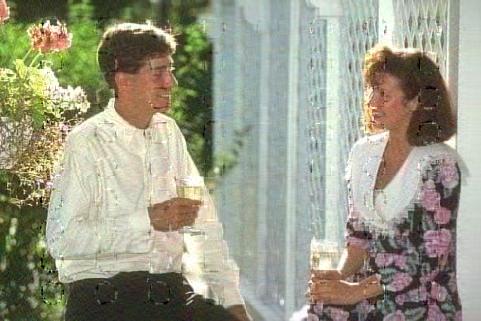}%
	\\
	%\text{{\footnotesize (b)}}
\end{minipage}%
\begin{minipage}[f]{0.123\linewidth}
	\centering
	\vspace{6.6pt}
	\includegraphics[width=2.15cm]{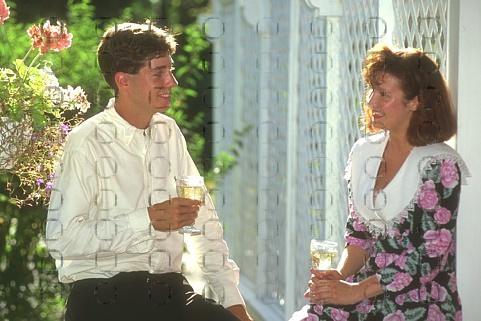}%
	\\
	%\text{{\footnotesize (c)}}
\end{minipage}%
\begin{minipage}[f]{0.123\linewidth}
	\centering
	\vspace{6.6pt}
	\includegraphics[width=2.15cm]{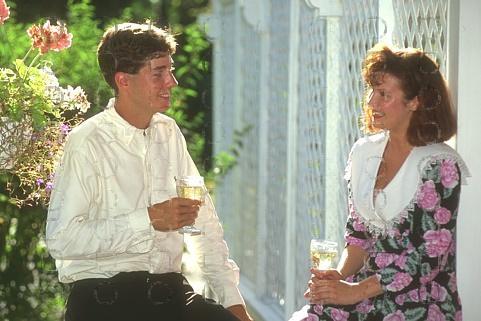}%
	\\
	%\text{{\footnotesize (d)}}
\end{minipage}%
\begin{minipage}[f]{0.123\linewidth}
	\centering
	\vspace{6.6pt}
	\includegraphics[width=2.15cm]{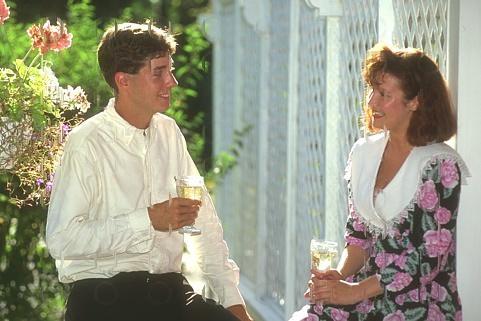}%
	\\
	%\text{{\footnotesize (e)}}
\end{minipage}%
\begin{minipage}[f]{0.123\linewidth}
	\centering
	\vspace{6.6pt}
	\includegraphics[width=2.15cm]{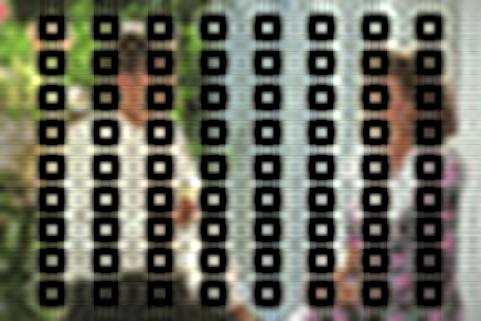}%
	\\
	%\text{{\footnotesize (f)}}
\end{minipage}%
\begin{minipage}[f]{0.123\linewidth}
	\centering
	\vspace{6.6pt}
	\includegraphics[width=2.15cm]{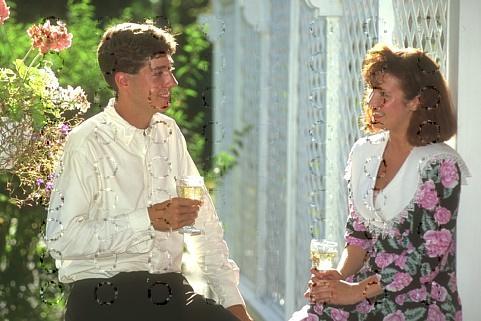}%
	\\
	%\text{{\footnotesize (g)}}
\end{minipage}%
\begin{minipage}[f]{0.123\linewidth}
	\centering
	\vspace{6.6pt}
	\includegraphics[width=2.15cm]{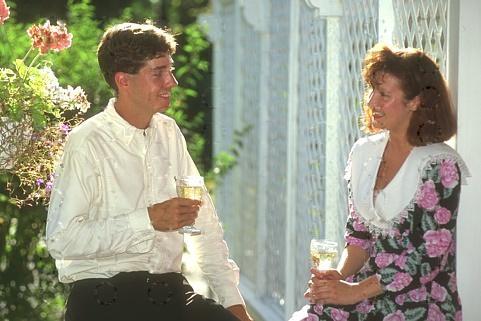}%
	\\
	%\text{{\footnotesize (h)}}
\end{minipage}%
\\
\vspace{1.5pt}%%%%%%%%%%%%%%%%%%%%%%%%%%%%%%%%%%%%%%%%%%%%%%%%%%%%
	\centering
\begin{minipage}[f]{0.123\linewidth}
	\centering
	\text{Alphabet Mask}
	\includegraphics[width=2.15cm]{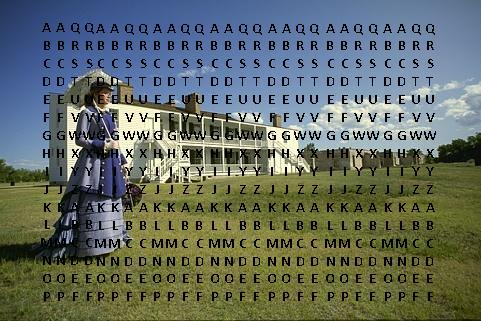}%
	\\
	\text{{\footnotesize (a)}}
\end{minipage}%
\begin{minipage}[f]{0.123\linewidth}
	\centering
	\vspace{6.6pt}
	\includegraphics[width=2.15cm]{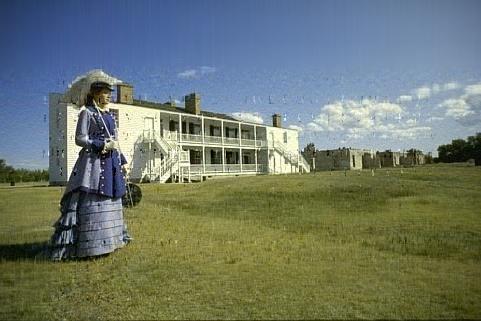}%
	\\
	\text{{\footnotesize (b)}}
\end{minipage}%
\begin{minipage}[f]{0.123\linewidth}
	\centering
	\vspace{6.6pt}
	\includegraphics[width=2.15cm]{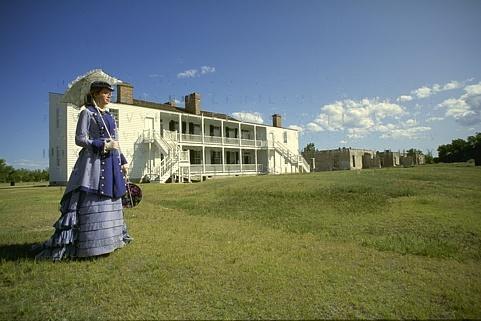}%
	\\
	\text{{\footnotesize (c)}}
\end{minipage}%
\begin{minipage}[f]{0.123\linewidth}
	\centering
	\vspace{6.6pt}
	\includegraphics[width=2.15cm]{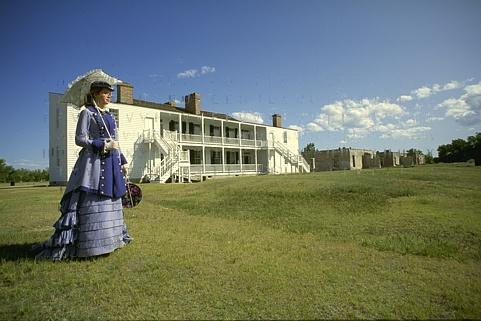}%
	\\
	\text{{\footnotesize (d)}}
\end{minipage}%
\begin{minipage}[f]{0.123\linewidth}
	\centering
	\vspace{6.6pt}
	\includegraphics[width=2.15cm]{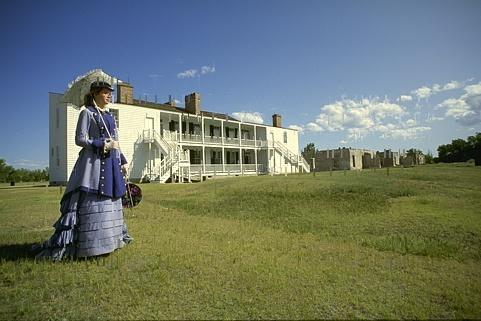}%
	\\
	\text{{\footnotesize (e)}}
\end{minipage}%
\begin{minipage}[f]{0.123\linewidth}
	\centering
	\vspace{6.6pt}
	\includegraphics[width=2.15cm]{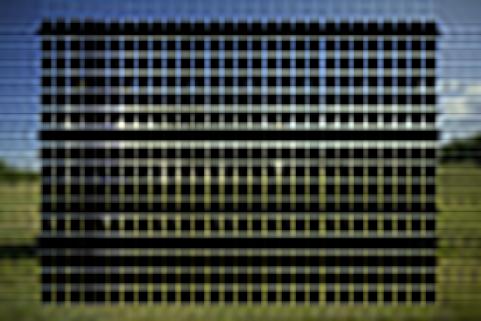}%
	\\
	\text{{\footnotesize (f)}}
\end{minipage}%
\begin{minipage}[f]{0.123\linewidth}
	\centering
	\vspace{6.6pt}
	\includegraphics[width=2.15cm]{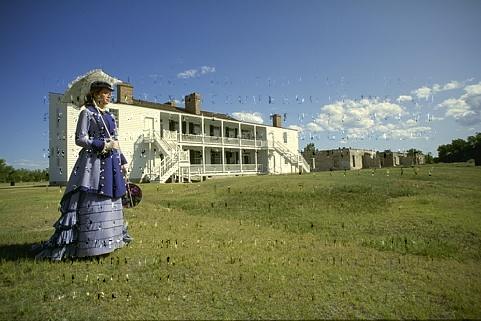}%
	\\
	\text{{\footnotesize (g)}}
\end{minipage}%
\begin{minipage}[f]{0.123\linewidth}
	\centering
	\vspace{6.6pt}
	\includegraphics[width=2.15cm]{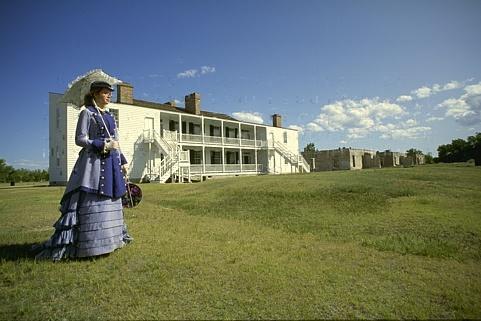}%
	\\
	\text{{\footnotesize (h)}}
\end{minipage}%
\\
\vspace{1.5pt}%%%%%%%%%%%%%%%%%%%%%%%%%%%%%%%%%%%%%%%%%%%%%%%%%%%%
	\caption{The observations and the recovered images via algorithms with different masks. Illustrations from left to right: (a) {Observations}, (b) {FBCP}, (c) {HaLRTC}, (d) {TNN}, (e) {SiLRTC-TT}, (f) {STDC}, (g) MF-TV, (h) HPMF.}
	\label{maskCI}
\end{figure*}

In Fig. \ref{maskCI}, we illuminate visual quality on the image recovery with masks. From the figure, except for STDC, all methods roughly remove different masks, and reconstruct image outlines. Nevertheless, the proposed HPMF recovers clearer image details and textures than compared algorithms. Concretely, low-rank prior based methods (FBCP, HaLRTC) shows limited details, especially when removing masks of Grid, Line, and Circle. This may tell that, masks along the tensor modes severely degrade the recovery accuracy of FBCP and HaLRTC. Other low-rank prior based methods (TNN, SiLRTC-TT) give sharp visual quality for all mask removals. For smoothness prior regularized methods (STDC, MF-TV), we see that, STDC fails to recover images corrupted by masks. MF-TV only presents clear results with masks of Text, Circle, and Alphabet. To sum up, our HPMF demonstrates the best visual quality to remove all masks on average in the experiments.

\section{Conclusion}
\label{sec-conclu}
In this article, we propose a novel matrix factorization model to hierarchically employ the low-rank prior, TV prior, and SC prior for tensor completion. The ADMM based algorithm is then utilized to solve the proposed model, and we also investigate the algorithm complexity and convergence. Experiments on various data sets demonstrate the superiority of our algorithm over several state-of-the-art algorithms. For tensor completion, the proposed  HPMF framework can be generalized by using different factorization methods and priors. Furthermore, we also expect this framework is applicable to other applications, e.g., the image denoising, the background subtraction, etc.

% or




\ifCLASSOPTIONcaptionsoff
  \newpage
\fi

\bibliographystyle{ieeetr}
\bibliography{HPMF_PrePrintRel_0410_v2.bbl}

%\input{HPMF_vFastRel_0311.bbl}
\end{document}